\pdfoutput=1
\documentclass{mnarticle}
\usepackage{amssymb}
\usepackage{graphicx}
\usepackage{xspace}
\usepackage{multirow}
\usepackage{mydef}
\usepackage[numbers]{natbib}
\let\cite=\citet
\usepackage{enumerate}
\usepackage{booktabs}
\usepackage{subcaption}
\usepackage{placeins}

\newsiamremark{remark}{Remark}

\begin{document}

\title{Maximum principle preserving approximation \\ of
  scalar conservation equations\\ with a fractional Laplacian}

\author{Andrea Bonito\thanks{Department of Mathematics, Texas A\&M University, College Station, TX 77843, USA (\email{bonito@tamu.edu}). The work of AB was partially supported by NSF grant DMS-2409807.}
\and Jean-Luc Guermond\thanks{Department of Mathematics, Texas A\&M University, College Station, TX 77843, USA (\email{guermond@tamu.edu}).}
\and Murtazo Nazarov\thanks{Department of Information Technology, Division of Scientific Computing, Uppsala University, Box 337, SE-751 05 Uppsala, Sweden (\email{murtazo.nazarov@uu.se}).
The work of MN was supported by the Swedish Research Council (VR) grant 2021-04620.}%
}

\headers{Monotone approximation of conservation equations with fractional Laplacian}{A. Bonito, J.-L. Guermond, and M. Nazarov}

\maketitle

\begin{abstract}
In this paper, we propose a maximum principle preserving continuous finite element method for scalar conservation laws augmented with a fractional diffusion operator of order $s\in (0,1)$. Conservation laws with fractional diffusion operators arise in a variety of physical models involving anomalous transport and nonlocal regularization, including geophysical flows, plasma physics, and nonlinear acoustics. The fractional Laplacian is approximated using a sinc-quadrature representation, yielding a resolvent-based formulation compatible with continuous finite elements and explicit strong stability preserving Runge--Kutta time integration. A first-order monotone scheme is constructed using graph viscosity and is shown to satisfy a discrete maximum principle under an explicit CFL condition. A high-order extension is obtained via entropy-consistent corrections and convex limiting, achieving second-order accuracy while preserving the maximum principle, and hence positivity, as well as mass conservation in the periodic setting. Numerical experiments, including applications to the Surface Quasi-Geostrophic equation, confirm the second-order accuracy and the robustness of the proposed method on unstructured meshes.

\end{abstract}



\section{Introduction}
\label{Sec:introduction}

Scalar conservation laws augmented with nonlocal diffusion operators arise in a wide range of physical models, including overdriven detonations in gases~\citep{Clavin_2002}, anomalous diffusion in crystal growth and semiconductor processes~\citep{Woyczynski_2001}, nonlinear acoustics~\citep{Sugimoto_1991}, geophysical flows~\citep{Constantin_1994,Ohkitani_1997}, and plasma transport~\citep{del_Castillo_2010}. In many of these systems, the underlying transport processes deviate substantially from classical Fickian diffusion, exhibiting long-range interactions or heavy-tailed jump processes that are more accurately captured by fractional diffusion operators. As a consequence, the fractional Laplacian has become a central tool for modeling anomalous transport and nonlocal regularization.

We consider scalar conservation laws with a fractional diffusion term of order $s\in(0,1)$, written as
\begin{equation}
\begin{cases}
  \partial_t u + \DIV\bef(u) + (-\LAP)^su = g,  &\quad (\bx,t) \in D\times\polR^+, \\
  u(\bx,t)= 0, &\quad (\bx,t) \in \partial D \times \polR^+,\\
  u(\bx,0)=u_0(\bx), &\quad \bx\in D,
  \label{PDE_system}
\end{cases}
\end{equation}
where $D\subset \polR^d$, $d=1,2,3$, is an open bounded domain, $u_0$ is an initial condition, $\bef:\polR\to\polR^d$ is the flux function, and $g(\bx,t)$ is a given scalar-valued source term. Here $(-\Delta)^s$ is the \emph{spectral} fractional Laplacian, see Section~\ref{Sec:Laps}.

Despite its wide applicability, the numerical approximation of \eqref{PDE_system} presents several challenges. The fractional Laplacian is intrinsically nonlocal, leading to dense discretization matrices and a high computational cost. In addition, solutions of problems involving fractional diffusion often develop steep gradients or exhibit limited regularity, complicating the design of numerical schemes that are accurate and robust.

A number of numerical methods have been proposed to approximate fractional diffusion operators, including finite differences by~\cite{Huang_2014}, spectral and Fourier-based approaches by~\cite{Constantin_2012,Constantin_1994}, spectral vanishing viscosity techniques by~\cite{Song_2017}, the Balakrishnan method by \cite{bonito2015numerical,bonito2017numerical,bonito2017sinc} and the extension method by \cite{nochetto2015pde}; we refer, for instance, to \citep{bonito2018numerical,d2020numerical} for surveys. However, when coupled with nonlinear hyperbolic fluxes, additional difficulties arise: numerical oscillations may violate the discrete maximum principle, positivity, or the invariant domain structure of the PDE. Ensuring monotonicity becomes even more delicate when the diffusion operator is nonlocal.

The primary goal of this work is to construct and analyze a space--time discretization for scalar conservation laws with a fractional Laplacian that \emph{preserves the maximum principle} at the discrete level. Our approach is based on continuous, piecewise-linear finite elements in space and strong-stability-preserving (SSP) explicit Runge--Kutta schemes in time. A key ingredient is a sinc-quadrature representation of the fractional Laplacian  \citep{KBSinc92} applied to the Balakrishnan-type integral identity
\begin{equation}\label{fractional_laplacian}
(-\Delta)^s u = c_s \int_{-\infty}^{\infty} e^{s y} (e^{y} -\Delta)^{-1} (-\Delta u)\, dy,
\end{equation}
which reduces the nonlocal operator to a weighted family of resolvent problems for the classical Laplacian. This representation enables a natural finite-element approximation that interacts cleanly with the underlying discrete structure of the method.
The representation \eqref{fractional_laplacian} was originally advocated by \cite{bonito2019numerical} for the integral fractional Laplacian and subsequently by \cite{bonito2020electroconvection} for the spectral fractional Laplacian.

To enforce the maximum principle, we construct a first-order monotone scheme obtained by combining the resolvent-based discretization of the fractional Laplacian with a graph viscosity mechanism, see, \eg \citep{Guermond_Nazarov_2014, guermond_popov_second_order_2018, Bonito_2021}. When the mesh satisfies the acute-angle condition, we prove that this low-order scheme satisfies a discrete maximum principle under an explicit CFL constraint that accounts for both the hyperbolic and fractional diffusive contributions.
This closes a gap left open in \cite{Bonito_2021}, where a discrete maximum principle was established for the inviscid SQG model only, while in the presence of fractional diffusion only an invariant-domain property was obtained under a more restrictive time-step condition \citep[Remark~3.3]{Bonito_2021}. The key observation is that, under the acute-angle condition, the discrete resolvents appearing in the sinc quadrature are monotone, so that the fractional diffusion term can be written as a convex combination of nodal values; see Lemma~\ref{Lem:min_max_Wl}, Lemma~\ref{l:local_form}, and Theorem~\ref{Thm:Min_Max_first_order}.
We also treat Dirichlet boundary conditions, which were not considered in \cite{Bonito_2021} and require a modification in the analysis.

Building upon this foundation, we design a high-order method obtained by adding entropy-consistent high-order corrections and applying the convex limiting framework of~\cite{Guermond_Nazarov_Popov_Tomas_2018}. The resulting algorithm achieves second-order accuracy in both space and time while preserving the maximum principle and, in the periodic setting, mass conservation. The limiter acts only where needed, allowing the high-order method to remain active in smooth regions without sacrificing non-oscillatory behavior near sharp layers.

The performance of the proposed method is illustrated through several benchmarks. We first verify the expected convergence rates for the fractional heat equation and for a manufactured Burgers' problem with fractional diffusion. We then apply the algorithm to the Surface Quasi-Geostrophic (SQG) equation, a widely studied geophysical model whose solutions may develop strong fronts. The method robustly captures the formation of sharp gradients while preserving the maximum principle, even on unstructured meshes.

The paper is organized as follows. In \S\ref{Sec:Laps} we review the sinc-quadrature approximation of the fractional Laplacian. The fully discrete scheme and the first-order invariant-domain preserving approximation are described in \S\ref{Sec:FE}, together with the CFL analysis. High-order extensions based on SSP time stepping and convex limiting are presented in \S\ref{Sec:Higher-order}. Numerical experiments are provided in \S\ref{sec:numerics}. Concluding remarks are given in \S\ref{sec:conclusions}.

\section{\texorpdfstring{The spectral fractional Laplacian $(-\LAP)^s$ and its approximations}{The spectral fractional Laplacian and its approximations}} \label{Sec:Laps}
We recall that $(-\Delta) : H^2(D) \cap H^1_0(D) \subset L^2(D) \rightarrow L^2(D)$ is an unbounded operator and has a compact inverse. Therefore, there is a countable collection of eigenpairs $\{ \lambda_{i} , e_{i}\}_{i \in \mathbb N} \subset \mathbb R^+ \times H^1_0(D)$ such that $\{ e_{i} \}_{i \in \mathbb N}$ constitutes an orthonormal basis of $L^2(D)$ as well as an orthogonal basis of $H^1_0(D)$. The spectral fractional Laplacian $(-\Delta)^s$, $s\in (0,1)$, is defined for smooth and compactly supported functions $u \in C^\infty_0(D)$ by
\begin{equation}\label{e:frac_lap_def}
(-\Delta)^s u := \sum_{i \in \mathbb N} \lambda_i^s u_i e_i, \qquad \textrm{with }u_i:= \int_{D} u e_i \quad \textrm{for }i\in \mathbb N.
\end{equation}
Furthermore, this definition extends to functions $u$ in the interpolation space
$H^s_0(D):= [L^2(D),H^1_0(D)]_s$ (complex interpolation method) by density.

In this section, we advocate for the equivalent representation \eqref{fractional_laplacian} of the spectral fractional Laplacian. In Section~\ref{s:heuristic} we give a heuristic justification of \eqref{fractional_laplacian} and refer to \cite{lunardi2018interpolation} for a proof using functional calculus.
Then, we approximate $(-\LAP)^s$ by introducing a quadrature
for \eqref{fractional_laplacian}. The key here is that, once
discretized, \eqref{fractional_laplacian} only involves the Laplacian
operator, which is then easy to handle using standard finite element
techniques.  The material in this section has been developed in
\citep{bonito2015numerical,bonito2017numerical,bonito2017sinc,bonito2019numerical,bonito2020electroconvection}. We also mention that the case of periodic boundary conditions, which are used in the numerical experiments of Section~\ref{sec:numerics}, can be treated similarly.
We do not pursue this option to avoid adding an additional layer of abstraction, but refer to \cite{Bonito_2021} and \cite{bonito2021approximation}.
The reader familiar with this material is invited to jump to \S\ref{Sec:FE}.

\subsection{Heuristics}\label{s:heuristic}
One way to understand the origin of the
expression~\eqref{fractional_laplacian} consists of realizing that the
following identity holds true for any $\lambda>0$ and any
$s\in(0,1)$:
\begin{equation}
  \lambda^{s} = \lambda^{s}  c_{s} \int_0^\infty \frac{z^{s-1}}{z+1}\diff z =
  c_{s} \int_{-\infty}^\infty  e^{s y} (e^y+\lambda)^{-1}\lambda \diff y, \label{lambda_identity}
\end{equation}
with $c_s^{-1} = \int_0^\infty \frac{z^{s-1}}{z+1}\diff z$ (use
the change of variables $\lambda z= e^y$).  Using the beta function
$B(x,y) = \int_0^\infty \frac{t^{x-1}}{(1+t)^{x+y}}\diff t=
\frac{\Gamma(x)\Gamma(y)}{\Gamma(x+y)}$,
defined for $\Re(x)>0$ and $\Re(y)>0$, we observe that
$c_s^{-1}=B(s,1-s)$. Then, Euler's reflection formula
$\Gamma(1-z)\Gamma(z)=\frac{\pi}{\sin(\pi z)}$, which holds for any $z\not\in\polZ$,
gives $c_s =\frac{\sin(\pi s)}{\pi}$. Using the basis $\{e_{i}\}_{i\in \mathbb N}$ made of eigenfunctions of $(-\Delta)$ introduced above and in view of \eqref{e:frac_lap_def}, we have
$$
(-\Delta)^s u = \sum_{i \in \mathbb N} u_i (-\Delta)^s e_i =  \sum_{i \in \mathbb N} u_i \lambda_i^s e_i
$$
so that \eqref{lambda_identity}, applied with $\lambda=\lambda_i$, gives \eqref{fractional_laplacian}.

\subsection{Generic quadrature}
Let $M, N$ be two positive integers. We set $\calL:=\intset{-M}{N}$,
with the notation
$\intset{a}{b}:=\{a,a+1,\ldots,b-1,b\}$
for any $a,b\in \polZ$ such that $a\le  b$.
Let $\{(\omega_l,y_l)\}_{l\in\calL}$ be the weights and nodes of a quadrature for $\int_{-\infty}^{\infty} f(t)\diff t$, then
we approximate $(-\LAP)^{s}u$ using \eqref{fractional_laplacian} as follows:
\begin{equation}
  (-\LAP)^{s}u \approx c_s \sum_{l\in\calL} \omega_l e^{s y_l} (e^{y_l}-\LAP)^{-1} (-\LAP u).
\end{equation}
In order to simplify this expression, we follow \cite{bonito2019numerical} and define
the mapping  $u\mapsto w_l(u)$ by setting $(e^{y_l}-\LAP)^{-1} (-\LAP u) = w_l(u)+u$, \ie
\begin{align*}
  -\LAP u &= (e^{y_l}-\LAP)(w_l(u)+u) =  (e^{y_l}-\LAP)w_l(u) +e^{y_l}u -\LAP u.
\end{align*}
In other words, the mapping $w_l: H^{-1}(\Dom)\to H^1_0(\Dom)$ is defined by setting
\begin{equation}
  (e^{y_l}-\LAP)w_l(u) = -e^{y_l} u, \label{def_wl}
\end{equation}
for any $u\in H^{-1}(\Dom)$. Hence,
$$
(-\Delta)^s u \approx c_s \sum_{l \in \mathcal L} \omega_l e^{s y_l} (u+w_l(u)).
$$
This leads to the following approximation of~\eqref{PDE_system}:
\begin{equation}
  \partial_t u + \DIV\bef(u) +\sum_{l\in\calL}  c_s\omega_l e^{s y_l} (w_l(u)+u) =g,\quad
  u_{|\front}=0, \quad u(\cdot,0)=u_0(\cdot).
  \label{approx_PDE_system}
\end{equation}
More details on the way the sinc quadrature is constructed and its accuracy
are given in the next section.

\subsection{Accuracy of the sinc quadrature}

We fix $k>0$, $M, N \in \mathbb N$ positive.
The sinc quadrature corresponds to the choice $\omega_l := k$ and $y_l:=l k$ for $l\in\calL:=\intset{-M}{N}$.

We now recall results about the sinc quadrature approximation
$$
k\sum_{l=-M}^N \left(e^{sy_l}(e^{y_l}-\Delta)^{-1}(-\Delta w), v \right)_{L^2(D)}
$$
of
$$
\int_{-\infty}^\infty \left(e^{sy}(e^y-\Delta)^{-1}(-\Delta w), v \right)_{L^2(D)} dy
$$
for $w \in H^{\delta}(D)$ with $\delta >s$ and $v \in H^s(D)$.
Here $(\cdot,\cdot)_{L^2(D)}$ denotes the $L^2$-scalar product and we let
$$
q(y;w,v) := \left(e^{sy}(e^y-\Delta)^{-1}(-\Delta w), v \right)_{L^2(D)}.
$$

The sinc quadrature formula takes advantage of (i) the analyticity of $z \mapsto q(z;w,v)$ on the strip
$$
D_{r} := \left\lbrace z \in \mathbb C \ : \  |\Im(z)| < r \right\rbrace,
$$
where $0< r < \pi$ is fixed from now on, and (ii) the decay property
$$
|q(z;w,v)| \leq C(r) \left\lbrace
\begin{array}{ll}
  e^{s\Re(z)} \|  w \|_{L^2(D)} \| v \|_{L^2(D)}, \qquad & \textrm{when } \Re (z) < 0, \\
  e^{(s-\delta)\Re(z)/2} \|  w \|_{H^\delta(D)} \| v \|_{H^s(D)}, \qquad & \textrm{when } \Re (z) \geq  0,
\end{array}
\right.
$$
where $C(r)$ is a constant only depending on $r$ which tends to infinity as $r \to \pi$.
From now on $C(r)$ denotes such a generic constant which may change from one relation to the next.

From (ii) we deduce that
$$
\int_{-r}^{r} | q(x+iy; w,v) | dy \leq C(r) \| w \|_{H^\delta(D)} \|v \|_{H^s(D)}
$$
and
$$
\int_{-\infty}^\infty \left(  | q(x+ir; w,v)| + |q(x-ir;w,v)| \right)dx \leq \frac{C(r)}{\min(\delta-s,s)}  \| w \|_{H^\delta(D)} \|v \|_{H^s(D)}.
$$
Whence, Theorem 2.20 in \cite{KBSinc92} applies to guarantee
$$
\left|\int_{-\infty}^\infty q(y;w,v)dy - k \sum_{j=-\infty}^\infty q(kj;w,v)\right| \leq C(r) e^{-2\pi r/k}  \| w \|_{H^\delta(D)} \|v \|_{H^s(D)}.
$$
In addition, the decay property (ii) yields estimates for the tails
$$
k \sum_{j= -\infty}^{-M-1} | q(kj;w,v)| \leq \frac{C(r)}{s} e^{-sMk}  \| w \|_{L^2(D)} \|v \|_{L^2(D)}
$$
and
$$
k \sum_{j=N+1}^{\infty} | q(kj;w,v)| \leq \frac{C(r)}{\delta-s} e^{(s-\delta)Nk/2}   \| w \|_{H^\delta(D)} \|v \|_{H^s(D)}.
$$
Gathering all these estimates, we find that upon selecting
\begin{equation}\label{choice_N_M}
M = \left\lceil \frac{2\pi r}{s k^2} \right\rceil \qquad \textrm{and} \qquad 
N=\left\lceil \frac{4\pi r}{(\delta-s)k^2} \right\rceil
\end{equation}
we arrive at
\begin{equation}\label{e:quad_error}
  \left| \int_{-\infty}^\infty  q(y;w,v) \, dy -
  k\sum_{j=-M}^N  q(kj;w,v) \right| \leq C(r,s,\delta) e^{-2\pi r/k}   \| w \|_{H^\delta(D)} \|v \|_{H^s(D)}.
\end{equation}

\section{Full approximation} \label{Sec:FE} In this section we describe
 the time and space approximation of~\eqref{approx_PDE_system}.

\subsection{Finite element setting}
Let $\famTh$ be a shape-regular sequence of subdivisions made of simplices (the case of quadrilaterals/hexahedra can be treated similarly). For simplicity, we assume that $\Dom$ is a polyhedron that is exactly matched by the subdivision $\mathcal T_h$, $h>0$.
We denote by $\wK$ the reference simplex and by $T_K:\wK \to K$ the generic affine map from the reference element to the physical element.

We approximate the solution of
\eqref{approx_PDE_system} with piecewise linear continuous Lagrange
finite elements generated from $\famTh$. The approximation space is
defined by
\begin{align} \label{eq:Xh}
  P(\calT_h) &=\{ v\in \calC^0(\Dom;\Real) \cap H^1_0(D) \st
  v_{|K}{\circ}T_K \in \mathbb P_1,\ \forall K\in \calT_h\},
\end{align}
where $\mathbb P_1$ is the set of polynomials of total degree no larger than 1. We denote by $\calV^+$ the index set of all the vertices of $\calT_h$, by $\calV^\partial\subset\calV^+$ the index set of the vertices located on $\partial\Dom$, and by $\calV:=\calV^+\setminus\calV^\partial$ the index set of the interior vertices. Let $\{\varphi_i\}_{i\in\calV^+}$ be the Lagrange basis functions associated with all the vertices. Then $\{\varphi_i\}_{i\in\calV}$ is a basis of $P(\calT_h)$, so that $\textrm{card}(\calV)=\dim P(\calT_h)$, and $\sum_{i\in\calV^+}\varphi_i=1$ (partition of unity). For any $i\in\calV$, we
let $\calI(i):=\{j\in \calV^+\st \varphi_i\varphi_j \not\equiv 0\}$ be the stencil of $i$; note that $\calI(i)$ contains the boundary vertices adjacent to $i$.
The coefficients of the functions in $P(\calT_h)$ are extended by zero on $\calV^\partial$: for $u_h=\sum_{j\in\calV}\sfU_j\varphi_j\in P(\calT_h)$, we set $\sfU_j:=0$ for $j\in\calV^\partial$, so that $u_h=\sum_{j\in\calV^+}\sfU_j\varphi_j$. This is the homogeneous Dirichlet condition in \eqref{PDE_system}. With this convention, the sums over $j\in\calI(i)$ below include the boundary vertices adjacent to $i$, where the coefficients vanish.

For any $i\in\calV$ and $j\in\calV^+$
we set $a_{ij}:=\int_\Dom \GRAD\varphi_i\SCAL\GRAD\varphi_j\diff x$. The partition of unity property implies $\sum_{j\in\calI(i)}a_{ij}=0$ for all $i\in\calV$.

We assume that each subdivision $\mathcal T_h$, $h>0$, is such that
\begin{equation}
  a_{ij}:=\int_\Dom \GRAD\varphi_i\SCAL\GRAD\varphi_j\diff x\le 0,\qquad
  \forall i\in\calV,\ \forall j\in \calI(i){\setminus}\{i\}. \label{stifness_matrix}
\end{equation}
This condition is known to be satisfied when $d=1$ and for simplices when $d=2$ or $d=3$ under the so-called acute angle
condition \cf \eg \cite{BrKKS:09}, \cite[\S5.2]{BrKoK:08},
\cite[Eq.~(2.5)]{XuZik:99}. The purpose of this paper is not to
relax this condition. Note that, since $\calI(i)$ contains the boundary vertices adjacent to $i$, \eqref{stifness_matrix} also applies to the pairs formed by an interior and a boundary vertex.

We henceforth denote by $\calM^L$ the lumped mass matrix, \ie the diagonal matrix with entries
\begin{equation}
  m_i:= \int_{\Dom} \varphi_i(\bx)\diff \bx,  \qquad i\in\calV.
\end{equation}
Notice that $\varphi_i(\bx)\ge 0$ for all $\bx\in \Dom$ and so $m_i>0$ for all
$i\in \calV$.

We also introduce the following two quantities which will play an
important role in the rest of the paper
\begin{equation}
  \bc_{ij} := \int_\Dom \varphi_i \GRAD\varphi_j \diff
  x,\qquad \bn_{ij}:= \frac{\bc_{ij}}{\|\bc_{ij}\|_{\ell^2}}\qquad
  i\in\calV,\ j\in\calV^+,
  \label{def_of_cij}
\end{equation}
where $\|\cdot\|_{\ell^2}$ stands for the Euclidean norm in $\Real^d$.
The partition of unity implies
$\sum_{j\in\calV^+} \bc_{ij} =\mathbf{0}$, \ie $\sum_{j\in\calI(i)}\bc_{ij}=\mathbf{0}$ for $i\in\calV$. 

\subsection{Approximation of \texorpdfstring{$w_l$}{wl}} We start by approximating the
mapping $w_l$ defined in \eqref{def_wl} for any $l\in \calL$.  Let
$u_h = \sum_{i\in\calV} \sfU_i \varphi_i$ be any function in
$P(\calT_h)$. We define the discrete mapping $w_h^l: P(\calT_h) \to P(\calT_h)$ by setting
$w_h^l(u_h) := \sum_{i\in\calV} \sfW_i^l(u_h) \varphi_i$, where the coefficients
$\{\sfW_i^l(u_h)\}_{i\in\calV}$ solve the following linear system:
\begin{equation}\label{e:Wli}
  m_i e^{y_l} \sfW_i^l(u_h) + \sum_{j\in\calI(i)} a_{ij} \sfW_j^l(u_h) =  -m_i e^{y_l} \sfU_i,\qquad i\in\calV, \qquad \textrm{and} \qquad \sfW_i^l(u_h) = 0, \qquad i \in \mathcal V^\partial.
\end{equation}
Recall that $a_{ij}$ are the coefficients of the stiffness matrix defined
in~\eqref{stifness_matrix} and $m_i$ those of the lumped mass matrix.
When there is no ambiguity, we write $\sfW_j^l$ instead of $\sfW_j^l(u_h)$.

We have the following maximum principle.
\begin{lemma}[Maximum principle] \label{Lem:min_max_Wl}
  Let $\calT_h$ be such that \eqref{stifness_matrix} holds.
  Let $u_h = \sum_{i\in \calV^+} \sfU_i \varphi_i \in P(\calT_h)$ with $\sfU_i = 0$ when $i \in \mathcal V^\partial$ and $w_h=\sum_{i\in\calV^+} \sfW_i^l(u_h) \varphi_i \in P(\calT_h)$.
  Then we have for any $l \in \calL$
  \begin{equation}
    -\max_{j\in\calV^+} \sfU_j \le \min_{j\in\calV^+}\sfW_j^l(u_h)
    \le\max_{j\in\calV^+} \sfW_j^l(u_h)\le -\min_{j\in\calV^+} \sfU_j.
  \end{equation}
\end{lemma}
\begin{proof}
  Using that $\sum_{j\in\calI(i)}a_{ij}=0$, which follows from the partition of unity property satisfied by $\{\varphi_i\}_{i\in \calV^+}$, we
  deduce
  \begin{align*}
    m_i e^{y_l} \sfW_i^l =  -m_i e^{y_l} \sfU_i
    +\sum_{j\in\calI(i){\setminus}\{i\}}  - a_{ij}(\sfW_j^l -\sfW_i^l),\qquad i\in\calV.
  \end{align*}
  Let $i\in \calV{^+}$ be the index where $\sfW_j^l$ is maximum over $j \in \calV{^+}$, \ie
  $\sfW_i^l=\max_{j\in\calV^+} \sfW_j^l$. If $i\in\calV^\partial$, then $\max_{j\in\calV^+}\sfW_j^l=0\le-\min_{j\in\calV^+}\sfU_j$, since $\sfU_j=0$ for $j\in\calV^\partial$. Instead, for $i\in\calV$, since $-a_{ij}\ge 0$
  for all $j\in\calI(i){\setminus}\{i\}$ thanks to  \eqref{stifness_matrix},
  we have
  \[
  \max_{j\in\calV^+} m_i \sfW_j^l=m_i \sfW_i^l \le -m_i \sfU_i.
  \]
  This means $\max_{j\in\calV^+} \sfW_j^l\le -\min_{j\in\calV^+} \sfU_{j}$.
  Similarly, let $i\in \calV^+$ be the index where $\sfW_j^l$ is minimum over $j\in \calV^+$,
  \ie $\sfW_i^l=\min_{j\in\calV^+} \sfW_j^l$. If $i\in\calV^\partial$, then $\min_{j\in\calV^+}\sfW_j^l=0\ge-\max_{j\in\calV^+}\sfU_j$. Instead, when $i\in\calV$ we have
  \[
  \min_{j\in\calV^+} m_i\sfW_j^l=m_i \sfW_i^l \ge -m_i \sfU_i.
  \]
  Altogether this implies $\min_{j\in\calV^+} \sfW_j^l\ge -\max_{j\in\calV^+} \sfU_{j}$.
\end{proof}

We will also use the following local representation of $\sfW_i^l+\sfU_i$, which involves the quantities $\Lambda_i := a_{ii}/m_i>0$, $i\in\calV$.

\begin{lemma}[Local form of the fractional term]\label{l:local_form}
Let $\calT_h$ be such that \eqref{stifness_matrix} holds.
Let $u_h=\sum_{i\in\calV}\sfU_i\varphi_i\in P(\calT_h)$ and $\sfW^l_j:=\sfW^l_j(u_h)$. Then, for all $i\in\calV$ and $l\in\calL$,
\begin{equation}\label{eq:local_form}
  -(\sfW_i^l+\sfU_i) = \sum_{j\in\calI(i)\setminus\{i\}} \theta_{ij}^l\,\big(-\sfW_j^l-\sfU_i\big),
  \qquad
  \theta_{ij}^l := \frac{-a_{ij}/m_i}{e^{y_l}+\Lambda_i}\ \ge 0,
\end{equation}
and $\sum_{j\in\calI(i)\setminus\{i\}}\theta_{ij}^l = \frac{\Lambda_i}{e^{y_l}+\Lambda_i}\in(0,1)$.
\end{lemma}

\begin{proof}
Dividing the $i$-th equation defining $\sfW^l$ \eqref{e:Wli} by $m_i e^{y_l}$ and using $\sum_{j\in\calI(i)}a_{ij}=0$, we obtain
\[
  \sfW_i^l+\sfU_i = -\frac{e^{-y_l}}{m_i}\sum_{j\in\calI(i)} a_{ij}\,\sfW_j^l
  = \frac{e^{-y_l}}{m_i}\sum_{j\in\calI(i)\setminus\{i\}} (-a_{ij})\big(\sfW_j^l-\sfW_i^l\big).
\]
We write $\sfW_j^l-\sfW_i^l=(\sfW_j^l+\sfU_i)-(\sfW_i^l+\sfU_i)$ and use
$\sum_{j\in\calI(i)\setminus\{i\}}(-a_{ij})=a_{ii}$ to get
\[
  \big(1+e^{-y_l}\Lambda_i\big)\big(\sfW_i^l+\sfU_i\big)
  =\frac{e^{-y_l}}{m_i}\sum_{j\in\calI(i)\setminus\{i\}}(-a_{ij})\big(\sfW_j^l+\sfU_i\big),
\]
which is \eqref{eq:local_form}. The weights $\theta_{ij}^l$ are nonnegative owing to \eqref{stifness_matrix},
and their sum is $\frac{a_{ii}/m_i}{e^{y_l}+\Lambda_i}$.
\end{proof}

Observe that if $i$ is an index corresponding to a vertex adjacent to one in $\partial\Dom$, the sum in \eqref{eq:local_form} contains the boundary vertices in $j\in\calI(i)\cap\calV^\partial$, for which $-\sfW_j^l=0$. The homogeneous boundary value thus enters in \eqref{eq:local_form} with the weights $\theta_{ij}^l$, $j\in\calI(i)\cap\calV^\partial$.

\subsection{Full approximation, low-order scheme}
\label{sec:low-order}
We now construct an approximation in time and space
of~\eqref{approx_PDE_system}.  The space approximation is performed using the finite element space $P(\calT_h)$ while the forward Euler method is used for the time discretization.

Let $\dt>0$ and set $t_n:= n \tau$, $n\in \polN$. Let $G^n:= \sum_{i\in \calV} \sfG_i^n \varphi_i$ be an approximation of $g$ at time $t=t_n$ and $u_h^n:=\sum_{i\in\calV} \sfU_i^n \varphi_i$ be a (possibly high-order) approximation of $u$ at time $t_n$.
We define the low-order approximation of $u$ at time $t_{n+1}$ by $u_h\upLnp=\sum_{i\in\calV} \sfU_i\upLnp \varphi_i$, where the coefficient $\sfU_i\upLnp$, $i \in \calV$, is given by the relation
\begin{equation}\label{eq:LO}
  \begin{aligned}
    m_i \frac{\sfU_i\upLnp - \sfU_i^{n}}{\dt} + \sum_{j\in \mathcal I(i)} \bc_{ij} \SCAL\bef(\sfU_j^n)
    &+ c_s m_i\sum_{l\in\calL} \omega_l e^{s y_l} (\sfW_i^{l,n}(u_h^n)+\sfU_i^{n}) \\
    &= m_i \sfG_i^n
    + \sum_{j\in\calI(i){\setminus}\{i\}} d_{ij}\upLn (\sfU_j^{n}-\sfU_i^{n}).
  \end{aligned}
\end{equation}
The artificial viscosity coefficients $d_{ij}\upLn$ originally proposed in \citep{Guermond_Popov_2016} are added to deal with the instabilities inherent to the approximations of advection systems. They are given for $i \in \calV$, $j\in \mathcal I(i) \setminus \{ i \}$ by
\begin{equation}\label{eq:dijL}
  d_{ij}\upLn:=
  \max\Big(
  \lambda_{\max}(\sfU_i^{n},\sfU_j^{n},\bn_{ij})\|\bc_{ij}\|_{\ell^2}, \
  \lambda_{\max}(\sfU_j^{n},\sfU_i^{n},\bn_{ji})\|\bc_{ji}\|_{\ell^2}
  \Big),
\end{equation}
where $\lambda_{\max}(\sfU_i^{n},\sfU_j^{n},\bn_{ij})$ is an upper bound of the maximum wave speed in the 1D Riemann problem
$\p_t v +\p_x (\bef(v)\bn_{ij}) =0$, with left and right data
$\sfU_i^n$ and $\sfU_j^n$ respectively. In view of the discussion in Section~79.2 in \cite{Ern_Guermond_2021}, for instance for convex fluxes $\bef$, we have
$$
\lambda_{\max}(\sfU_i^{n},\sfU_j^{n},\bn_{ij}) =
\begin{cases}
\frac{| \bef(\sfU_i^{n})\SCAL\bn_{ij} - \bef(\sfU_j^{n})\SCAL\bn_{ij}|}{| \sfU_i^{n} - \sfU_j^{n}|} & \qquad \textrm{if} \quad \sfU_i^{n} > \sfU_j^{n},\\
\max \big(| \bef'(\sfU_i^{n})\SCAL\bn_{ij}|, |\bef'(\sfU_j^{n})\SCAL\bn_{ij}| \big) & \qquad \textrm{otherwise}.
\end{cases}
$$
For any $\bef$ (not necessarily convex), this choice for $\lambda_{\max}$ guarantees that
  \[
  \overline \sfU_{ij}^{n} := \frac12(\sfU_j^{n}+\sfU_i^{n})
  -(\bef(\sfU_j^n)-\bef(\sfU_{i}^n)) \SCAL\frac{\bc_{ij}}{2d_{ij}\upLn}
  \]
  satisfies
  \begin{equation}\label{e:interm_state}
  \overline \sfU_{ij}^{n}\in [\min(\sfU_i^{n},\sfU_j^{n}),\max(\sfU_i^{n},\sfU_j^{n})];
  \end{equation}
  see \citep[Lemma~79.18, Remark~79.19]{Ern_Guermond_2021}.

When $i=j$, we set $d_{ii}\upLn:= - \sum_{j\in\calI(i){\setminus}\{i\}} d_{ij}\upLn$, which implies the following important properties of the artificial viscosity coefficients:
$$
d_{ij}\upLn=d_{ji}\upLn,\quad  d_{ij}\upLn \ge 0 \qquad \forall i,j, \ i\neq j, \qquad \textrm{and} \qquad  \sum_{j\in\calI(i)} d_{ij}\upLn = 0 \quad \forall i.
$$

Recall that from \eqref{lambda_identity}, $\Lambda_i^s = c_s\int_{-\infty}^{\infty} e^{sy}\,(e^y+\Lambda_i)^{-1}\, \Lambda_i \diff y$, where $\Lambda_i=a_{ii}/m_i$ is introduced before Lemma~\ref{l:local_form}. 
This motivates the notation for $i\in\calV$
\begin{equation}\label{eq:Lambda_q}
  \Lambda_{i,\textrm{sinc}}^s  := c_s\sum_{l\in\calL}\omega_l\, e^{s y_l}(e^{y_l}+\Lambda_i)^{-1}\Lambda_i,
\end{equation}
to denote the sinc approximation of $\Lambda_i^s$.

The next result guarantees a local discrete maximum principle satisfied by the coefficients $\sfU_i\upLnp$ under the time-step restriction
\begin{equation}\label{eq:CFL}
\dt :=  \frac{\textrm{CFL}}{\max_{i\in\calV}\Big( \frac{1}{m_i}\sum_{j\in\calI(i){\setminus}\{i\}}d_{ij}\upLn
  + \Lambda_{i,\textrm{sinc}}^s\Big)},
\end{equation}
where $\textrm{CFL}>0$ is a constant.
We postpone to Section~\ref{s:CFL} a discussion of the relation between \eqref{eq:CFL} and the standard CFL condition for the transport equation.

We are now in a position to prove a discrete maximum property. 
Before proceeding, we recall the convention that $\sfU_j^n=\sfW_j^{l,n}=0$ whenever $j \in \mathcal V^\partial$ indexes a boundary vertex. These boundary values enter the discrete maximum estimates stated in the following result.

\begin{theorem}[Discrete Maximum Principle]\label{Thm:Min_Max_first_order}
  Let $u_h^n=\sum_{i\in \cal V} \sfU_i^n \varphi_i$ be an approximation at time $t_n$ and $\sfU_i\upLnp$ be as in \eqref{eq:LO}. Assume that the time-step $\dt$ satisfies \eqref{eq:CFL} for some CFL number $0 < \textrm{CFL} \le \frac12$.
  Then for any $i\in \calV$ we have
  \begin{equation}\label{eq:UL_max}
    \sfU_i^{n,\min} + \dt\sfG_i^n\le \sfU_i\upLnp
    \le \sfU_i^{n,\max} + \dt\sfG_i^n,
  \end{equation}
  where
\begin{equation}\label{Umax-Umin}
\begin{aligned}
  \sfU_i^{n,\min}&:=\min\Big(\min_{j\in\calI(i)}\sfU_j^{n},\ \min_{l\in\calL}\min_{j\in\calI(i)\setminus\{i\}}\big({-}\sfW_j^{l,n}\big)\Big),\qquad \textrm{and} \\
  \sfU_i^{n,\max}&:=\max\Big(\max_{j\in\calI(i)}\sfU_j^{n},\ \max_{l\in\calL}\max_{j\in\calI(i)\setminus\{i\}}\big({-}\sfW_j^{l,n}\big)\Big).
\end{aligned}
\end{equation}
Here $\sfW_{j}^{l,n}:=\sfW_{j}^l(u_h^n)$.
\end{theorem}

\begin{proof}
  Let $i\in \calV$. Using $\sum_{j\in\calI(i)}\bc_{ij}=\bzero$, we have
  \begin{align*}
    \sfU_i\upLnp ={}& \sfU_i^{n} +\frac{\dt}{m_i}\sum_{j\in\calI(i){\setminus}\{i\}} \Big(-\bc_{ij}\SCAL
    (\bef(\sfU_j^n)-\bef(\sfU_i^n)) + d_{ij}\upLn (\sfU_j^{n}-\sfU_i^{n})\Big) \\
    & - \dt c_s\sum_{l\in\calL} \omega_l e^{s y_l} (\sfW_i^{l,n}+\sfU_i^{n}) + \dt\sfG_i^n.
  \end{align*}
  In view of \eqref{e:interm_state}, we find that
  \begin{equation}\label{expr_UiLn}
  \begin{split}
    \sfU_i\upLnp ={}& \sfU_i^{n}(1-2\frac{\dt}{m_i}\sum_{j\in\calI(i){\setminus}\{i\}}d_{ij}\upLn)
    +\sum_{j\in\calI(i){\setminus}\{i\}} \frac{\dt}{m_i}2d_{ij}\upLn\overline \sfU_{ij}^{n} \\
    & - \dt c_s\sum_{l\in\calL} \omega_l e^{s y_l} (\sfW_i^{l,n}+\sfU_i^{n})
    + \dt\sfG_i^n.
    \end{split}
  \end{equation}

Owing to Lemma~\ref{l:local_form} and the definition \eqref{eq:Lambda_q} of $\Lambda_{i,\textrm{sinc}}^s$, the contribution of the fractional Laplacian in \eqref{expr_UiLn} can be rewritten as
\[
 -\dt c_s\sum_{l\in\calL}\omega_l e^{sy_l}\big(\sfW_i^{l,n}+\sfU_i^n\big)
 = \sum_{l\in\calL}\sum_{j\in\calI(i)\setminus\{i\}}\dt\,c_s\,\omega_l e^{sy_l}\theta_{ij}^l\big(-\sfW_j^{l,n}\big) - \dt\, \Lambda_{i,\textrm{sinc}}^s\,\sfU_i^n .
\]
Inserting this identity into \eqref{expr_UiLn}, we obtain
\begin{multline*}
  \sfU_i\upLnp-\dt\sfG_i^n
  = \Big(1-\frac{2\dt}{m_i}\sum_{j}d_{ij}\upLn-\dt\, \Lambda_{i,\textrm{sinc}}^s\Big)\sfU_i^n
  + \frac{2\dt}{m_i}\sum_{j}d_{ij}\upLn\,\overline\sfU_{ij}^n 
  + \sum_{l\in\calL}\sum_{j}\dt\,c_s\,\omega_l\, e^{sy_l}\,\theta_{ij}^l\,\big(-\sfW_j^{l,n}\big),
\end{multline*}
where all the sums over $j$ run over $\calI(i)\setminus\{i\}$. From Lemma~\ref{l:local_form}, we find that
$$
\sum_{j \in \calI(i)\setminus\{i\}} \theta_{ij}^l = \frac{\Lambda_i}{e^{y_l}+\Lambda_i}
$$
and thus the coefficients in front of $\sfU_i^n$, $\overline\sfU_{ij}^n$, and $-\sfW_j^{l,n}$ sum to one.
They are nonnegative because \eqref{eq:CFL} with $\textrm{CFL}\le\frac12$ implies
$\frac{2\dt}{m_i}\sum_{j\ne i}d_{ij}\upLn+\dt\, \Lambda_{i,\textrm{sinc}}^s\le 2\,\textrm{CFL}\le 1$.
Hence $\sfU_i\upLnp-\dt\sfG_i^n$ is a convex combination of $\sfU_i^n$, $\overline\sfU_{ij}^n$, $-\sfW_j^{l,n}$. Since $\overline\sfU_{ij}^n\in[\min(\sfU_i^n,\sfU_j^n),\max(\sfU_i^n,\sfU_j^n)]$ according to \eqref{e:interm_state}, we obtain \eqref{eq:UL_max} for $\sfU_i^{n,\min}, \sfU_i^{n,\max}$ given by \eqref{Umax-Umin}.
\end{proof}

\begin{remark}[Local and global bounds]\label{r:local_global}
The limits \eqref{Umax-Umin} involve, in addition to the values of $u_h^n$ on the stencil $\calI(i)$, the neighboring nodal values $-\sfW_{j}^{l,n}$, $j\in\calI(i)\setminus\{i\}$, $l\in\calL$, which are weighted averages of $u_h^n$ over the domain. This is unavoidable since $(-\LAP)^s$ is nonlocal and a maximum principle involving only the stencil $\calI(i)$ cannot hold in general. According to Lemma~\ref{Lem:min_max_Wl}, we have $\min_{j\in\calV^+}\sfU_j^n\le\sfU_i^{n,\min}\le\sfU_i^{n,\max}\le\max_{j\in\calV^+}\sfU_j^n$, so that \eqref{eq:UL_max} implies the global maximum principle $\min_{j\in\calV^+}\sfU_j^n+\dt\sfG_i^n\le\sfU_i\upLnp\le\max_{j\in\calV^+}\sfU_j^n+\dt\sfG_i^n$; recall that $\sfU_j^n=0$ for $j\in\calV^\partial$. The relations \eqref{Umax-Umin} are available at a negligible extra cost since $\sfW_{j}^{l,n}$ are needed anyway to evaluate \eqref{eq:LO}.
\end{remark}

\subsection{Estimation of the time-step restriction}\label{s:CFL}
For the local discrete maximum principle to hold (Theorem~\ref{Thm:Min_Max_first_order}), the time-step $\dt$ is required to satisfy \eqref{eq:CFL}.
Note that the first term in the denominator is the advection contribution, and the second term, $\Lambda_{i,\textrm{sinc}}^s$, corresponds to the fractional diffusion contribution. The following lemma shows that the sinc quadrature $\Lambda_{i,\textrm{sinc}}^s$ of $\Lambda_i^s$, see \eqref{eq:Lambda_q}, behaves like $\Lambda_i^s$ independently of the sinc quadrature parameters $M$ and $N$.

\begin{lemma}[Bound on $\Lambda_{i,\textrm{sinc}}^s$]\label{l:q_bound}
For all $M,N\in\polN$ and $k>0$, if the sinc quadrature weight and points are given by $\omega_l=k$ and $y_l=lk$, $l=-M,...,N$, then 
\[
  \Lambda_{i,\textrm{sinc}}^s \le c_s\Big(\frac{1}{s(1-s)}+2k\Big)\Lambda_i^s \le \frac{4+2k}{\pi}\,\Lambda_i^s, \qquad i\in\calV.
\]
\end{lemma}

\begin{proof}
We use $\frac{\Lambda}{e^y+\Lambda}\le\min(1,\Lambda e^{-y})$. Let $L\in\polZ$ be such that
$e^{Lk}\le\Lambda_i<e^{(L+1)k}$. Then
\[
  \frac{\Lambda_{i,\textrm{sinc}}^s}{c_s}\le \sum_{l\le L} k\, e^{slk} + \Lambda_i\sum_{l\ge L+1} k\, e^{(s-1)lk}
  = \frac{k\,e^{sLk}}{1-e^{-sk}} + \frac{k\,\Lambda_i\, e^{(s-1)(L+1)k}}{1-e^{-(1-s)k}}
  \le \Lambda_i^s\Big(\frac{k}{1-e^{-sk}}+\frac{k}{1-e^{-(1-s)k}}\Big),
\]
where we used $e^{sLk} \leq \Lambda_i^s$ and $e^{(s-1)(L+1)k} \leq \Lambda_i^{s-1}$.
The first bound follows from $\frac{x}{1-e^{-x}}\le 1+x$ for $x>0$, applied with $x=sk$ and $x=(1-s)k$. The second follows from
$c_s=\frac{\sin(\pi s)}{\pi}\le\frac1\pi$ and $c_s=\frac{\sin(\pi s)}{\pi}\le \frac{4s(1-s)}{\pi}$ for $s\in(0,1)$.
\end{proof}

A few comments are in order. 
Since $\Lambda_i=a_{ii}/m_i\le C h^{-2}$, where $h$ denotes the smallest element diameter (the mesh size for quasi-uniform subdivisions $\mathcal T_h$) and $C$ only depends on the subdivision shape-regularity constant,  Lemma~\ref{l:q_bound} implies that $\Lambda_{i,\textrm{sinc}}^s\le C h^{-2s}$ for a constant $C$ uniform in $h$, $M$ and $N$. Hence the CFL condition \eqref{eq:CFL} becomes  
\begin{equation}\label{e:dt_cond}
\dt\lesssim\min(h,h^{2s}), 
\end{equation}
which is the expected stability condition. In particular, for $s\le\frac12$ it is compatible with the usual CFL condition for transport equations. 

We also observe that the condition $\dt\, \Lambda_{i,\textrm{sinc}}^s\lesssim1$ included in \eqref{eq:CFL} with $\Lambda_{i,\textrm{sinc}}^s\approx(a_{ii}/m_i)^s$ is the fractional analogue of the classical condition $\dt\, a_{ii}/m_i\le1$ for the explicit approximation of the heat equation.

Finally, because the stability condition \eqref{e:dt_cond} is uniform in the sinc parameters $M$ and $N$, their values only affect the accuracy. Therefore, they can be chosen as in \eqref{choice_N_M} to balance the different quadrature errors without affecting the stability of the discrete scheme. 

\section{High-order accuracy}
\label{Sec:Higher-order}

In this section, we show how to increase the accuracy in time and space
of the method presented in Section~\ref{sec:low-order}.

\subsection{Time approximation}
Increasing the time accuracy of the method described in the above
section while keeping the maximum principle can be done by using
so-called strong stability preserving (SSP) time discretization
methods.  The key idea is to achieve higher-order accuracy in time by
making convex combinations of forward Euler steps. More precisely, each
time step of a SSP method is decomposed into substeps that are all
forward Euler steps; the final update is constructed as a convex
combination of the intermediate solutions. We refer the reader to
\citep{Ferracina_Spijker_2005,Higueras_2005,Gottlieb_Ketcheson_Shu_2009}
for detailed reviews on SSP methods.

Without claiming originality, let us illustrate the SSP concept. Let
us consider a finite-dimensional vector space $E$, a subset
$A\subset E$ and a (possibly nonlinear) operator
$L:[0,T]\CROSS A\longrightarrow E$.  We are interested in
approximating in time the following system of ordinary
differential equations $\partial_t u =L(t,u)$ with appropriate
initial condition.  We assume that this problem is well-posed (for
instance $L$ is continuous w.r.t. $t$ and Lipschitz w.r.t. $u$).  We
further assume that there exists a convex subset $B\subset A$ and
$\dt_{\max}>0$ such that
\begin{equation}
  v+ \dt L(t,v) \in  B,\qquad \forall v\in B,\ \forall t\in [0,T],\
  \forall \dt \in [0,\dt_{\max}]. \label{Forward_Euler_assumption}
\end{equation}

Let us now consider an $S$-stage
explicit Runge--Kutta method, and let us assume that this method has
the following $(\alpha-\beta)$ representation: There are real
coefficients $\alpha_{ik}$, $\beta_{ik}$ with $k\in\intset{0}{i-1}$
and $i\in\intset{1}{S}$ such that $u^{n+1}$ is obtained by first
setting $w^{(0)} = u^n$, then computing
\begin{align}
  w^{(i)} &= \sum_{k=0}^{i-1} \alpha_{ik}w^{(k)} +
  \beta_{ik}\dt L(t^n+\gamma_{k}\dt,w^{(k)}),\quad i\in\intset{1}{S},\label{Shu_Osher_Representation_kth}
\end{align}
and finally setting $u^{n+1} = w^{(S)}$; here $t^n+\gamma_k\dt$ is the time at which stage $k$ is evaluated.  We further assume that
$\alpha_{ik}\ge0$, $\beta_{ik}\ge0$, $\sum_{k=0}^{i-1}\alpha_{ik}=1$, and that $\beta_{ik}=0$ if $\alpha_{ik}=0$, for all $k\in\intset{0}{i-1}$,
$i\in\intset{1}{S}$.  Not every $S$-stage explicit Runge--Kutta
method admits an $(\alpha-\beta)$ representation with nonnegative coefficients. Any Runge--Kutta
method that admits an $(\alpha-\beta)$ representation as defined above
is said to be SSP for the reason stated in the following theorem.
\begin{theorem}[Shu--Osher] \label{Th:Shu_Osher} Let
  $c_\textup{os}:=\inf_{\{\alpha_{ik}\not=0,\ 1\le k+1\le i\le
    S\}}\alpha_{ik}\beta_{ik}^{-1}$.
  Let $B\subset A$ be convex.  Let $u^n\in B$ and let $u^{n+1}$ be
  given by \eqref{Shu_Osher_Representation_kth}.  Assume that
  $\dt \le c_\textup{os} \dt_{\max}$, then $u^{n+1}\in B$.
\end{theorem}

For instance, the SSPRK$(2,2)$ method can be represented as follows:
\begin{equation}
  \begin{aligned}
    w^{(1)}&= u^n\! + \dt L(t^n,u^n), \\
    w^{(2)}&= w^{(1)}\! + \dt L(t^{n+1},w^{(1)}),\\
    u^{n+1}&= \tfrac12 u^n\! + \tfrac12 w^{(2)}.
  \end{aligned}
\end{equation}
The SSPRK$(3,3)$ method can be represented as follows:
\begin{equation}
  \begin{aligned}
    w^{(1)} &= u^n + \dt L(t^n,u^n),
    &&z^{(1)} = w^{(1)} + \dt L(t^n+\dt,w^{(1)}), \\
    w^{(2)} &= \frac34 u^n + \frac14 z^{(1)},\quad
    &&z^{(2)} = w^{(2)} + \dt L(t^n+\tfrac12\dt,w^{(2)}), \\
    u^{n+1} &= \frac13 u^n + \frac23 z^{(2)}.
  \end{aligned}
\end{equation}

\begin{remark}[Literature] We remark that Theorem~\ref{Th:Shu_Osher} has been
  established in a slightly different form in
  \citet[Prop.~2.1]{Shu_Osher1988} without explicitly invoking convexity.
  The variant of this theorem presented here relies only on convexity.
\end{remark}

\subsection{\texorpdfstring{Second-order approximation in space}{Second-order approximation in space}}
In Section~\ref{sec:low-order}, we introduced a low-order method to approximate \eqref{PDE_system} while ensuring that the maximum principle is preserved. The low-order scheme \eqref{eq:LO} serves as a crucial component in the limiting technique discussed in the subsequent section. Our objective is to construct a method that is at least second-order accurate and that also preserves the maximum principle. To achieve this goal, we require two additional components: a high-order entropy consistent method (that may violate the maximum principle due to oscillations in its solutions), and a reliable limiting algorithm that uses both the low- and high-order solutions to construct a higher-order maximum principle preserving approximation.

We fix $(\eta,\bq)$ to be an entropy pair for \eqref{PDE_system}, i.e., $\eta \in C^1(\mathbb R)$ is convex and $\bq \in C^1(\mathbb R;\mathbb R^d)$ satisfies $\bq'(u) = \eta'(u) \bef'(u)$.
The associated entropy residual is defined as
\[
R(u) :=
\big(
\partial_t u + \DIV\bef(u) + (-\Delta)^su -g
\big) \eta'(u) = \big(
\partial_t u + (-\Delta)^su -g
\big) \eta'(u) + \nabla \cdot \bq(u)
\]
and vanishes when $u$ is a smooth solution to \eqref{PDE_system}.

In order to derive a discrete counterpart of the above entropy residual, we follow the strategy adopted in  \cite{guermond2011entropy} and \cite{Guermond_Nazarov_Popov_Tomas_2018}.
Given $u_h^n = \sum_{i \in \calV} \sfU_i^n \varphi_i$, an approximation of $u$ at time $t_n$, we let
$\sfU_i\upGnp$ be defined as the low order approximation \eqref{eq:LO} but without stabilization terms, i.e., it satisfies
\begin{equation}\label{eq:GAL}
  \begin{aligned}
    m_i\frac{\sfU_i\upGnp - \sfU_i^{n}}{\dt} = -\sum_{j\in \mathcal I(i)} \bc_{ij} \SCAL\bef(\sfU_j^n)
    -c_s m_i\sum_{l\in\calL} \omega_l e^{s y_l} (\sfW_i^{l,n}(u_h^n)+\sfU_i^{n}) + m_i \sfG_i^n.
  \end{aligned}
\end{equation}
Then, the discrete entropy residual reads $R_h(u):= \sum_{i\in \calV}R_i^n \varphi_i$, where
\begin{equation}\label{Ri}
  \begin{aligned}
    R_i^{n}
    :=&\
    \Big(
    \frac{\sfU_i\upGnp - \sfU_i^{n}}{\dt}
    \Big) \eta'(\sfU^n_i) + \frac{1}{m_{i}}\sum_{j\in \mathcal I(i)} \bc_{ij} \SCAL\bq(\sfU_j^n)
    \\
    &\qquad +
    \Big(
    c_s \sum_{l\in\calL} \omega_l e^{s y_l} (\sfW_i^{l,n}(u_h^n)+\sfU_i^{n}) -  \sfG_i^n
    \Big) \eta'(\sfU^n_i).
  \end{aligned}
\end{equation}

After multiplying \eqref{eq:GAL} by $\eta'(\sfU_i^n)$, we notice that the coefficient $R_i^n$ in \eqref{Ri} reduces to
\begin{equation}\label{eq:entropy_commutator}
  \begin{aligned}
    R_i^{n}
    =&\
    - \frac{1}{m_{i}}\sum_{j\in \mathcal I(i)} \bc_{ij} \SCAL\bef(\sfU_j^n) \eta'(\sfU^n_i)
    + \frac{1}{m_{i}}\sum_{j\in \mathcal I(i)} \bc_{ij} \SCAL\bq(\sfU_j^n) \\
    =:&\
    R_{1,i}^{n} + R_{2,i}^{n}.
  \end{aligned}
\end{equation}

Next, we propose a nodal residual indicator for each node $i\in\calV$:
\begin{equation*}
\alpha_i^n := \frac{|R_{1,i}^n + R_{2,i}^n|}{|R_{1,i}^n| + |R_{2,i}^n| + \epsilon^R_\calV},
\end{equation*}
where $\epsilon^R_\calV:= 10^{-8}\max_{j \in \calV}\big(|R_{1,j}^n| + |R_{2,j}^n| \big)$ is a small positive constant introduced to avoid division by zero.
On the one hand, in regions where the solution is smooth, the entropy residual is small, \ie $R_{1,i}^n$ and $R_{2,i}^n$ nearly cancel, and so the residual indicator is small. On the other hand, where the solution is discontinuous, the entropy residual is large and the residual indicator is close to one.

We are now in position to define the high-order graph viscosity coefficient and the high-order approximation.
We set
\begin{equation}\label{eq:dijH}
  d_{ij}\upHn := d_{ij}\upLn \max(\alpha_i^n, \alpha_j^n) \quad \textrm{when } i\not = j, \quad d_{ii}\upHn:= -\sum_{j\in\calI(i)\setminus\{i\}} d_{ij}\upHn.
\end{equation}
Since $\alpha_j^n$ is only defined for $j\in\calV$, we set $\alpha_j^n:=0$ for $j\in\calV^\partial$, so that only the indicator of the interior vertex $i$ enters $d_{ij}\upHn$ when $j$ is a boundary vertex.
Note that, similarly to the low-order graph viscosity, we have $d_{ij}\upHn=d_{ji}\upHn$, $d_{ij}\upHn \ge 0$ for all $i \neq j$, and $\sum_{j\in\calI(i)} d_{ij}\upHn = 0$.

Given an approximation $u_h^n=\sum_{i\in\calV} \sfU_i^n \varphi_i$ of the solution $u$ to \eqref{approx_PDE_system} at time $t_n$, the high-order approximation $u_h\upHnp=\sum_{i\in\calV} \sfU_i\upHnp \varphi_i$ of $u$ at time $t_{n+1}$ is defined by the relations
\begin{equation}\label{eq:HO}
  \begin{aligned}
    \sum_{j \in \mathcal I(i)} m_{ij} \frac{\sfU_{j}\upHnp - \sfU_{j}^{n}}{\dt} + \sum_{j\in \mathcal I(i)} \bc_{ij} \SCAL\bef(\sfU_j^n)
    &+ c_s m_i\sum_{l\in\calL} \omega_l e^{s y_l} (\sfW_i^{l,n}(u_h^n)+\sfU_i^{n}) \\
    &= m_i \sfG_i^n
    + \sum_{j\in\calI(i){\setminus}\{i\}} d_{ij}\upHn (\sfU_j^{n}-\sfU_i^{n}).
  \end{aligned}
\end{equation}
Compared to the low-order approximation \eqref{eq:LO}, where a lumped mass matrix is used, here $m_{ij}:=\int_D \varphi_j(\bx) \varphi_i(\bx) \diff \bx$ are the coefficients of the consistent mass matrix.

\subsection{Relation between high- and low-order solutions}
To create a limiter matrix, it is necessary to establish the relationship between the high-order solution and the low-order solution while keeping in mind that the low-order solution maintains the maximum principle whereas the high-order solution may not.

Let us write the first-order and the high-order schemes as
\begin{align*}
  m_i\sfU_i\upLnp ={}& m_i\sfU_i^{n} +\dt\sum_{j\in\calI(i){\setminus}\{i\}} \big(-\bc_{ij}\SCAL
  \bef(\sfU_j^n) + d_{ij}\upLn (\sfU_j^{n}-\sfU_i^{n}) \big) \\
  & - \dt m_ic_s\sum_{l\in\calL} \omega_l e^{s y_l} (\sfW_i^{l,n}+\sfU_i^{n}) + \dt m_i\sfG_i^n,
\end{align*}
and
\begin{align*}
  \sum_{j\in \mathcal I(i)}m_{ij}\sfU_j\upHnp ={}& \sum_{j\in \mathcal I(i)}m_{ij}\sfU_j^{n} +\dt\sum_{j\in\calI(i){\setminus}\{i\}} \big(-\bc_{ij}\SCAL
  \bef(\sfU_j^n) + d_{ij}\upHn (\sfU_j^{n}-\sfU_i^{n})\big) \\
  & - \dt m_i c_s\sum_{l\in\calL} \omega_l e^{s y_l} (\sfW_i^{l,n}+\sfU_i^{n}) + \dt m_i \sfG_i^n.
\end{align*}
By subtracting the low-order scheme from the high-order scheme, we obtain:
\begin{align*}
  \sum_{j\in \mathcal I(i)}m_{ij}\sfU_j\upHnp
  ={}&
  m_i\sfU_i\upLnp
  +
  \sum_{j\in \mathcal I(i)}m_{ij}\sfU_j^{n} -  m_{i}\sfU_i^{n}
    +\dt\sum_{j\in\calI(i){\setminus}\{i\}} (d_{ij}\upHn - d_{ij}\upLn) (\sfU_j^{n}-\sfU_i^{n}).
\end{align*}
By observing that $m_i = \sum_{j\in \mathcal I(i)}m_{ij}$ we add $m_i \sfU_i\upHnp -\sum_{j\in\calI(i)}m_{ij}\sfU_j\upHnp$ to both sides of the last equality and obtain
\begin{align*}
  m_i\sfU_i\upHnp
  = &\
  m_i\sfU_i\upLnp
  + \sum_{j\in \mathcal I(i)}m_{ij}(\sfU_i\upHnp - \sfU_j\upHnp)
  + \sum_{j\in \mathcal I(i)}m_{ij}(\sfU_j^n - \sfU_i^n) \\
  &\ \qquad
  +\dt\sum_{j\in\calI(i){\setminus}\{i\}} (d_{ij}\upHn - d_{ij}\upLn) (\sfU_j^{n}-\sfU_i^{n}),
\end{align*}
which we write in the form
\begin{equation}\label{eq:H_L}
   m_{i}\sfU_i\upHnp =   m_i\sfU_i\upLnp + \sum_{j\in \mathcal I(i)} \bsfA_{ij}^{n+1},
\end{equation}
where for $i\in \calV$ and $j\in \mathcal I(i)$
\begin{align*}
  \bsfA_{ij}^{n+1} := (m_i \delta_{ij} - m_{ij}) \Big( \sfU_j\upHnp - \sfU_i\upHnp - (\sfU_j^{n}-\sfU_i^{n}) \Big) 
  + \dt (d_{ij}\upHn - d_{ij}\upLn) (\sfU_j^{n}-\sfU_i^{n}).
\end{align*}
The equation \eqref{eq:H_L} expresses the relation between the high-order and low-order solutions. 

We observe that the skew-symmetry property $\bsfA^{n+1}_{ij}=-\bsfA^{n+1}_{ji}$ for $i,j\in\calV$ implies
\begin{equation}\label{e:mass_high_low}
\sum_{i\in\calV}m_i\sfU_i\upHnp=\sum_{i\in\calV}m_i\sfU_i\upLnp+\sum_{i\in\calV}\sum_{j\in\calI(i)\cap\calV^\partial}\bsfA_{ij}^{n+1},
\end{equation}
since the contributions of the pairs of interior vertices cancel. In the periodic setting, $\calV^\partial=\emptyset$ and both schemes conserve the mass, \ie $\sum_{i\in\calV}m_i\sfU_i\upHnp=\sum_{i\in\calV}m_i\sfU_i\upLnp=\sum_{i\in\calV}m_i\sfU_i^n$. With homogeneous Dirichlet conditions, the mass is not conserved at the continuous level, since $\int_\Dom(-\LAP)^s u\diff\bx\neq0$ in general, and neither scheme conserves it: both exchange mass with the boundary through the pairs $(i,j)$ with $j\in\calV^\partial$. Relation \eqref{e:mass_high_low} shows that the high-order correction only redistributes mass among the interior vertices, except along these pairs, where it modifies the boundary flux of the low-order scheme by $\sum_{i\in\calV}\sum_{j\in\calI(i)\cap\calV^\partial}\bsfA_{ij}^{n+1}$.

\subsection{Convex limiting}

The objective of convex limiting, see \eg \citep[Sec.~4]{Guermond_Nazarov_Popov_Tomas_2018}, is to determine a limiter matrix $(\ell_{ij})_{i\in\calV,\,j\in\calI(i)}$ such that
$\ell_{ij} \in [0,1]$ for $i\in \calV$, $j\in \calI(i)$, and $\ell_{ij}=\ell_{ji}$ for $i \in \calV$, $j\in \calI(i)\cap\calV$ (symmetry is only required between interior degrees of freedom), and
\begin{equation}\label{eq:limited}
  m_{i}\sfU_i^{n+1} := m_i\sfU_i\upLnp + \sum_{j\in \mathcal I(i)} \ell_{ij}\bsfA_{ij}^{n+1},
\end{equation}
satisfies $\Psi_i^\pm(\sfU_i^{n+1}) \ge 0$, where
$$
\Psi_i^{-}(w):= \sfU_i^{n,\max} + \dt \sfG_i^n - w \qquad \textrm{and} \qquad \Psi_i^{+}(w):= w - \sfU_i^{n,\min} - \dt \sfG_i^n.
$$
Recall that $\sfU_i^{n,\max}$ and $\sfU_i^{n,\min}$ are defined in \eqref{Umax-Umin}.

It should be noted that when $\ell_{ij}\equiv 1$ in \eqref{eq:limited}, we obtain the high-order solution $\sfU_i\upHnp$, which may potentially violate the maximum principle. In contrast, setting $\ell_{ij}\equiv 0$ corresponds to the low-order solution $\sfU_i\upLnp$, which is known to satisfy the maximum principle. Additionally, similarly to \eqref{e:mass_high_low} but using in addition the symmetry assumption $\ell_{ij}=\ell_{ji}$ for $i \in\calV$ and $j\in \calI(i) \cap \calV$, we have
\[
\sum_{i\in\calV}m_i\sfU_i^{n+1}=\sum_{i\in\calV}m_i\sfU_i\upLnp+\sum_{i\in\calV}\sum_{j\in\calI(i)\cap\calV^\partial}\ell_{ij}\bsfA_{ij}^{n+1}.
\]

We now propose a limiter matrix and prove in Theorem~\ref{t:ho_max_p} that it leads to a discrete maximum principle.
Given $u_h^n=\sum_{i\in \calV} \sfU_i^n \varphi_i$ an approximation at time $t_n$ and $u_h\upLnp=\sum_{i\in \calV} \sfU_i\upLnp \varphi_i$ a low-order approximation at time $t_{n+1}$ satisfying the local discrete maximum principle \eqref{eq:UL_max}, for $i\in \calV$, we set $\ell_{ij}=0$ when $j\not \in \mathcal I(i)$ and
\begin{equation}\label{eq:lij}
  \ell_{ij}:= \min(\ell_j^i, \ell_i^j)\quad\textrm{if } j\in\calI(i)\cap\calV,\qquad \ell_{ij}:=\ell_j^i\quad\textrm{if } j\in\calI(i)\cap\calV^\partial,
\end{equation}
with
\begin{equation*}
  \ell_j^i :=
  \begin{cases}
    \frac{|\sfU_i^{n,\min} + \dt \sfG_i^n - \sfU_i\upLnp|}{|\bsfP_{ij}^{n+1}|} \!\!&  \mbox{if } \sfU_i\upLnp + \bsfP_{ij}^{n+1} < \sfU_i^{n, \min} + \dt \sfG_i^n, \\
    \frac{|\sfU_i^{n,\max} + \dt \sfG_i^n - \sfU_i\upLnp|}{|\bsfP_{ij}^{n+1}|} \!\!& \mbox{if } \sfU_i^{n, \max} + \dt \sfG_i^n < \sfU_i\upLnp + \bsfP_{ij}^{n+1},\\
    1 \!\!& \mbox{otherwise },
    \end{cases}
\end{equation*}
where $\bsfP_{ij}^{n+1} := \frac{1}{\lambda_j(i) m_i} \bsfA_{ij}^{n+1}$ with $\lambda_j(i):= \frac{1}{\textrm{card}(\calI(i)) - 1}$ for $i \in \calV$ and $j\in\calI(i)\setminus \{i\}$.
Note that for $j\in\calI(i)\cap\calV^\partial$, $\ell_i^j$ is not defined and the flux between $i$ and $j$ is then limited using only $\ell_j^i$.
Also, when $\sfU_i\upLnp + \bsfP_{ij}^{n+1} < \sfU_i^{n, \min} + \dt \sfG_i^n$, because the low-order approximation is assumed to satisfy the local discrete maximum principle \eqref{eq:UL_max}, we deduce that
\begin{equation}\label{e:Psign}
\bsfP_{ij}^{n+1} < \sfU_i^{n, \min} + \dt \sfG_i^n - \sfU_i\upLnp \leq 0.
\end{equation}
Therefore, dividing by $|\bsfP_{ij}^{n+1}|$ in the definition of $\ell_j^i$ is justified, and
$$
0\leq \ell_j^i=\frac{|\sfU_i^{n,\min} + \dt \sfG_i^n - \sfU_i\upLnp|}{|\bsfP_{ij}^{n+1}|}\leq 1.
$$
A similar argument shows that for the case $\sfU_i^{n, \max} + \dt \sfG_i^n < \sfU_i\upLnp + \bsfP_{ij}^{n+1}$, we have $\bsfP_{ij}^{n+1}>0$ and
$$
0\leq \ell_j^i=\frac{|\sfU_i^{n,\max} + \dt \sfG_i^n - \sfU_i\upLnp|}{|\bsfP_{ij}^{n+1}|} \leq 1.
$$
In practice, we avoid numerical instabilities by replacing $|\bsfP_{ij}^{n+1}|$ by $|\bsfP_{ij}^{n+1}|+ 10^{-14}\max_{i\in\calV}|\sfU_i^{n}|$ in the denominator.

\begin{theorem}\label{t:ho_max_p}
Let $u_h^n=\sum_{i\in \cal V} \sfU_i^n \varphi_i$ be an approximation at time $t_n$, $\sfU_i\upLnp$ the low-order approximation given by \eqref{eq:LO} and $\sfU_i^{n+1}$ be the convex limiting solution defined by \eqref{eq:limited} with the limiting matrix given in \eqref{eq:lij}.
 Assume that the time-step $\dt$ satisfies \eqref{eq:CFL} for some $0<\textrm{CFL}\le\frac12$.
 Then we have
   \begin{equation*}
    \sfU_i^{n,\min} + \dt\sfG_i^n\le \sfU_i^{n+1}
    \le \sfU_i^{n,\max} + \dt\sfG_i^n,
  \end{equation*}
for every $i\in \calV$.
\end{theorem}
\begin{proof}
We fix $i\in \cal V$ and remark that the stated bounds are equivalent to showing that
$\Psi_i^\pm(\sfU_i^{n+1}) \ge 0$. We only prove the result for $\Psi_i^+$ since the proof for $\Psi_i^-$ is analogous.

We start by noting that \eqref{eq:limited} can be rewritten as a convex combination
\begin{equation}\label{eq:conv}
  \sfU_i^{n+1} = \sum_{j\in\calI(i){\setminus}\{i\}} \lambda_j(i) \Big( \sfU_i\upLnp + \ell_{ij}\bsfP_{ij}^{n+1} \Big).
\end{equation}
Furthermore, since the upper level set $L(\Psi_i^+) := \{ r \in \polR \, |\,  \Psi_i^+(r)\ge 0 \}$ is convex, the  function $\Psi_i^+ : \polR \mapsto \polR$ is quasi-concave. Therefore, if
$
\Psi_i^{+}(\sfU_i\upLnp + \ell_{ij}\bsfP_{ij}^{n+1}) \ge 0
$
for all $j\in\calI(i){\setminus}\{i\}$, then the convex-combination representation \eqref{eq:conv} of $\sfU_i^{n+1}$ implies
$
\sfU_i^{n+1} \in L(\Psi_i^+),
$
i.e., $\Psi_i^+(\sfU_i^{n+1}) \ge 0$ as desired.

It remains to show that
$
  \Psi_i^{+}(\sfU_i\upLnp + \ell_{ij}\bsfP_{ij}^{n+1}) \ge 0,
$
which is performed by distinguishing three cases.\\

\noindent\fbox{\scriptsize 1} When $\sfU_i^{n, \min} \!+\! \dt \sfG_i^n  \le \sfU_i\upLnp \!+\! \bsfP_{ij}^{n+1} \le \sfU_i^{n, \max} \!+\! \dt \sfG_i^n$, we have
\[
  \Psi_i^{+}(\sfU_i\upLnp + \bsfP_{ij}^{n+1}) \ge 0.
\]
In addition, the maximum principle property \eqref{eq:UL_max} satisfied by the low-order coefficient $\sfU_i\upLnp$ reads
\[
  \Psi_i^{+}(\sfU_i\upLnp) \ge 0.
\]
Therefore, since $0\leq \ell_{ij} \leq 1$,
\[
	\sfU_i\upLnp + \ell_{ij}\bsfP_{ij}^{n+1} =
	(1-\ell_{ij})\sfU_i\upLnp + \ell_{ij} (\sfU_i\upLnp + \bsfP_{ij}^{n+1}),
\]
is a convex combination in $L(\Psi_i^+)$ and thus
\[
  \Psi_i^{+}(\sfU_i\upLnp + \ell_{ij}\bsfP_{ij}^{n+1}) \ge 0.
\]

\noindent\fbox{\scriptsize 2} When $\sfU_i\upLnp \!+\! \bsfP_{ij}^{n+1} < \sfU_i^{n, \min} \!+\! \dt \sfG_i^n$,
we recall that $\bsfP_{ij}^{n+1} < 0$ according to \eqref{e:Psign}.
Whence, from the definition of $\ell_{j}^{i}$, we observe that
\begin{align*}
  \Psi_i^{+}(\sfU_i\upLnp + \ell_j^i\bsfP_{ij}^{n+1})
&=
  \Psi_i^{+}\Bigg(
  \sfU_i\upLnp + \frac{\sfU_i\upLnp - \sfU_i^{n,\min} - \dt \sfG_i^n}{-\bsfP_{ij}^{n+1}}\,  \bsfP_{ij}^{n+1}
  \Bigg) \\
&=
  \Psi_i^{+}(
  \sfU_i^{n,\min} + \dt \sfG_i^n) = 0.
\end{align*}

\noindent\fbox{\scriptsize 3} When $\sfU_i^{n, \max} + \dt \sfG_i^n < \sfU_i\upLnp + \bsfP_{ij}^{n+1}$,
we have this time that $\bsfP_{ij}^{n+1} > 0$ and from the definition of $\ell_j^i$ we deduce that
\begin{align*}
  \Psi_i^{+}(\sfU_i\upLnp + \ell_j^i\bsfP_{ij}^{n+1})
&=
  \Psi_i^{+}\Bigg(
  \sfU_i\upLnp + \frac{\sfU_i^{n,\max} + \dt \sfG_i^n - \sfU_i\upLnp}{\bsfP_{ij}^{n+1}}\,  \bsfP_{ij}^{n+1}
  \Bigg) \\
&=
  \Psi_i^{+}(
  \sfU_i^{n,\max} + \dt \sfG_i^n) \ge 0.
\end{align*}
In cases 2 and 3, thanks to $0\leq \ell_{ij} \leq \ell_j^i$ and $0 \le \ell_j^i<1$ (if $\ell_j^i=0$, then $\ell_{ij}=0$ and there is nothing to prove),
\[
 \sfU_i\upLnp + \ell_{ij}\bsfP_{ij}^{n+1}
 =
 \Big(1-\frac{\ell_{ij}}{\ell_j^i} \Big) \sfU_i\upLnp +
 \frac{\ell_{ij}}{\ell_j^i} \Big( \sfU_i\upLnp + \ell_j^i\bsfP_{ij}^{n+1} \Big),
\]
is a convex combination between $\sfU_i\upLnp$ and $\sfU_i\upLnp + \ell_j^i\bsfP_{ij}^{n+1}$.
It remains to invoke again that the low-order approximation satisfies the discrete maximum principle, i.e., $\Psi_i^{+}(\sfU_i\upLnp) \geq 0$ to deduce that
$
  \Psi_i^{+}(\sfU_i\upLnp + \ell_{ij}\bsfP_{ij}^{n+1}) \ge 0.
$
\end{proof}

\section{Numerical examples}
\label{sec:numerics}
In this section, we test the performance of the proposed algorithm  on several benchmarks. The first benchmark is designed to test the fractional diffusion term while the second benchmark incorporates the Burgers-type hyperbolic term $\bef(u)=(\frac 1 2 u^2,\frac 1 2 u^2)^\intercal$. In both cases, exact solutions can be derived upon appropriately choosing $g$.

Finally, we apply the proposed method to the Surface Quasi-Geostrophic (SQG) equation. The SQG model is well known and widely used to study large-scale atmospheric and oceanic fluid dynamics. Common numerical approaches for solving SQG equations include spectral methods \citep{Constantin_1994,Constantin_2012}. A finite element method for the (mostly) inviscid SQG equation was first proposed in \citep{Bonito_2021}. By solving a standard benchmark problem, we validate that the proposed method is accurate and preserves the maximum principle for several values of the dissipation coefficient. 

The implementation uses the C++ interface of the FEniCS 
project \citep{fenics}. We use the third-order strong stability preserving explicit
Runge--Kutta method SSPRK$(3,3)$ of Section~\ref{Sec:Higher-order}.

\subsection{Fractional heat equation}\label{s:heat}
To validate our implementation of the fractional Laplacian, we consider the fractional heat equation
\begin{equation}
\begin{cases}
\begin{aligned}
\partial_t u(\bx,t) + \kappa (-\Delta)^s u(\bx,t) = 0, &\quad (\bx,t) \in D \times (0,T], \\
u(\bx, 0)
= u_0(\bx):=\sin x_2 \cos x_1,&  \quad \bx \in D,
\end{aligned}
\end{cases}
\end{equation}
with $D=[0, 2\pi]^2$ and periodic boundary conditions in both the $x_1$- and $x_2$-directions.
The time step is set according to \eqref{eq:CFL} with $\textrm{CFL}=0.4$.

In this setting, the exact solution is given by
\begin{align*}
u(\bx, t) := e^{-2^s \kappa t} \sin x_2 \cos x_1.
\end{align*}

We consider two values for the fractional power: $s = \frac{1}{4}$ and $s = \frac{1}{2}$. The results of the numerical computations are presented in Table~\ref{table:heat} for both the limited and unlimited cases. The unlimited case corresponds to the scheme \eqref{eq:HO}, while the limited solution is obtained using the convex limiting  \eqref{eq:limited}. The sequence of unstructured meshes is generated using the mesh generator \texttt{Gmsh}~\citep{gmsh_2009} and is constructed to satisfy the acute angle condition.

We observe second-order convergence in both the $L^1$- and $L^2$-norms for the limited and unlimited simulations (the rates are computed with respect to the mesh size $h\propto(\#\textrm{dofs})^{-1/2}$). Furthermore, the errors are nearly identical for both approaches on fine meshes, suggesting that the nonlinear viscosity approximation \eqref{eq:HO} essentially satisfies the bounds \eqref{eq:UL_max}, so that the limiter is rarely active in this regime.

Although not reported here, we observe that the convergence rate is not optimal in the $L^{\infty}$-norm for this set of unstructured meshes. This is in contrast to the optimal second-order rate observed in \citep[Sec.~4.2]{Bonito_2021}  when using uniform meshes.

\begin{table}[t]
  \caption{Fractional heat equation. Convergence history in the $L^1$- and $L^2$-norms for simulations using limited and unlimited approximations at $T = 0.5$ and $\kappa = 0.001$, with two different fractional diffusion powers $s = \frac{1}{4}$ and $s = \frac{1}{2}$.}
  \centering
  Unlimited solution \vspace{0.05in}\\
  \begin{tabular}{ccccc|cccc}
    \toprule
    & \multicolumn{4}{c}{$s=\frac12$}
    & \multicolumn{4}{c}{$s=\frac14$} \vspace{0.05in}\\
    \toprule
    \# dofs  & $L^1$ & rate    & $L^2$     & rate   & $L^1$ & rate    & $L^2$     & rate \\
     2700 &   4.89E-02 &   --   &   9.78E-03 &   --    &   5.16E-02 &   --   &   1.04E-02 &   --  \\
     5849 &   2.03E-02 &   2.28 &   4.03E-03 &   2.29  &   2.08E-02 &   2.35 &   4.16E-03 &   2.38\\
    13303 &   8.43E-03 &   2.13 &   1.71E-03 &   2.09  &   8.82E-03 &   2.09 &   1.79E-03 &   2.06\\
    29663 &   3.66E-03 &   2.08 &   7.85E-04 &   1.94  &   3.77E-03 &   2.12 &   8.11E-04 &   1.97\\
    66041 &   1.57E-03 &   2.12 &   3.56E-04 &   1.98  &   1.40E-03 &   2.48 &   3.37E-04 &   2.19\\
   148301 &   6.74E-04 &   2.09 &   1.39E-04 &   2.32  &   2.90E-04 &   3.89 &   8.54E-05 &   3.39\\
   333089 &   2.91E-04 &   2.08 &   6.08E-05 &   2.04  &   1.30E-04 &   1.98 &   3.84E-05 &   1.97\\
    \bottomrule
  \end{tabular}
  \\
  \vspace{0.05in}
  Limited solution
  \vspace{0.05in}\\
  \begin{tabular}{ccccc|cccc}
    \toprule
    & \multicolumn{4}{c}{$s=\frac12$}
    & \multicolumn{4}{c}{$s=\frac14$} \vspace{0.05in}\\
    \toprule
    \# dofs  & $L^1$ & rate    & $L^2$     & rate   & $L^1$ & rate    & $L^2$     & rate \\
     2700 &   4.90E-02 &   --   &   9.81E-03 &   --   &   5.23E-02 &   --   &   1.07E-02 &   --   \\
     5849 &   2.04E-02 &   2.27 &   4.07E-03 &   2.28 &   2.09E-02 &   2.38 &   4.22E-03 &   2.41 \\
    13303 &   8.45E-03 &   2.14 &   1.72E-03 &   2.10 &   8.93E-03 &   2.07 &   1.83E-03 &   2.03 \\
    29663 &   3.68E-03 &   2.07 &   7.98E-04 &   1.91 &   3.88E-03 &   2.08 &   8.71E-04 &   1.86 \\
    66041 &   1.57E-03 &   2.13 &   3.60E-04 &   1.99 &   1.42E-03 &   2.52 &   3.55E-04 &   2.24 \\
   148301 &   6.75E-04 &   2.09 &   1.40E-04 &   2.34 &   3.04E-04 &   3.80 &   9.91E-05 &   3.15 \\
   333089 &   2.91E-04 &   2.08 &   6.11E-05 &   2.05 &   1.33E-04 &   2.05 &   4.25E-05 &   2.09 \\
    \bottomrule
  \end{tabular}

  \label{table:heat}
\end{table}

\subsection{\texorpdfstring{Burgers' equation with fractional Laplacian}{Burgers' equation with fractional Laplacian}}
We consider the following two-dimensional fractional diffusion Burgers' equation:
\begin{equation}\label{eq:burger}
\begin{cases}
\begin{aligned}
  \partial_t u + \DIV \Big(\frac{u^2}{2}, \frac{u^2}{2} \Big)^\intercal + \kappa(-\Delta)^s u
&= g, \quad (\bx, t) \in D \times (0, T], \\
u(\bx, 0)
&= u_0(\bx), \quad \bx \in D,
\end{aligned}
\end{cases}
\end{equation}
where the domain is $D = [0, 2\pi]^2$ with periodic boundary conditions, $s = \frac{1}{2}$, and $T = 0.5$. The time step is set according to \eqref{eq:CFL} with $\textrm{CFL}=0.4$. In addition, we let
\begin{align*}
g(\bx, t) := \frac{1}{2}
e^{-2^{s+1} \kappa t} \Big(
- \sin (2 x_1) \sin^2 x_2
+ \cos^2 x_1 \sin (2 x_2)
\Big),
\end{align*}
corresponding to the manufactured solution
$
u(\bx, t) := e^{-2^s \kappa t} \sin x_2 \cos x_1.
$

The convergence history in two different norms is reported in Table~\ref{table:burger} for two diffusion coefficients: $\kappa = 0.01$ and $\kappa = 0.001$. The same sequence of unstructured meshes as in Section~\ref{s:heat} is used.
We observe second-order convergence rates in both norms, as reported in the table.

\begin{table}[t]
  \caption{Fractional diffusion Burgers' equation. Convergence history in the $L^1$- and $L^2$-norms for simulations using limited and unlimited approximations at $T = 0.5$ and $s=\frac12$, with two different fractional diffusion coefficients $\kappa = 0.01$ and $\kappa = 0.001$.}

  \centering
  Unlimited solution \vspace{0.05in}\\
  \begin{tabular}{ccccc|cccc}
    \toprule
    & \multicolumn{4}{c}{$\kappa = 0.01$}
    & \multicolumn{4}{c}{$\kappa = 0.001$} \vspace{0.05in}\\
    \toprule
    \# dofs  & $L^1$ & rate    & $L^2$     & rate   & $L^1$ & rate    & $L^2$     & rate \\
      2700 &   6.01E-02 &   --   &   1.21E-02 &   --   &   6.15E-02 &   --   &   1.24E-02 &   --   \\
      5849 &   2.70E-02 &   2.07 &   5.48E-03 &   2.05 &   2.84E-02 &   2.00 &   5.85E-03 &   1.95 \\
     13303 &   1.18E-02 &   2.01 &   2.42E-03 &   1.99 &   1.26E-02 &   1.97 &   2.67E-03 &   1.91 \\
     29663 &   5.07E-03 &   2.11 &   1.08E-03 &   2.01 &   5.77E-03 &   1.96 &   1.31E-03 &   1.78 \\
     66041 &   2.08E-03 &   2.22 &   4.16E-04 &   2.39 &   2.36E-03 &   2.24 &   5.01E-04 &   2.40 \\
    148301 &   9.32E-04 &   1.99 &   2.03E-04 &   1.77 &   1.14E-03 &   1.79 &   3.12E-04 &   1.17 \\
    333089 &   3.87E-04 &   2.17 &   8.26E-05 &   2.23 &   5.00E-04 &   2.04 &   1.47E-04 &   1.86 \\
    \bottomrule
  \end{tabular}
  \\
  \vspace{0.05in}
  Limited solution
  \vspace{0.05in}\\
  \begin{tabular}{ccccc|cccc}
    \toprule
    & \multicolumn{4}{c}{$\kappa = 0.01$}
    & \multicolumn{4}{c}{$\kappa = 0.001$} \vspace{0.05in}\\
    \toprule
    \# dofs  & $L^1$ & rate    & $L^2$     & rate   & $L^1$ & rate    & $L^2$     & rate \\
      2700 &   6.01E-02 &   --   &   1.21E-02 &   --   &   6.16E-02 &   --   &   1.25E-02 &   --   \\
      5849 &   2.71E-02 &   2.07 &   5.51E-03 &   2.04 &   2.85E-02 &   1.99 &   5.89E-03 &   1.94 \\
     13303 &   1.18E-02 &   2.02 &   2.42E-03 &   2.00 &   1.27E-02 &   1.98 &   2.68E-03 &   1.92 \\
     29663 &   5.08E-03 &   2.11 &   1.09E-03 &   2.00 &   5.79E-03 &   1.95 &   1.31E-03 &   1.78 \\
     66041 &   2.08E-03 &   2.23 &   4.16E-04 &   2.40 &   2.36E-03 &   2.24 &   5.03E-04 &   2.40 \\
    148301 &   9.32E-04 &   1.99 &   2.03E-04 &   1.77 &   1.14E-03 &   1.79 &   3.13E-04 &   1.17 \\
    333089 &   3.87E-04 &   2.17 &   8.26E-05 &   2.23 &   5.00E-04 &   2.04 &   1.48E-04 &   1.86 \\
    \bottomrule
  \end{tabular}
  \label{table:burger}
\end{table}

\subsection{\texorpdfstring{Surface Quasi-Geostrophic model}{Surface Quasi-Geostrophic model}}
Consider the following two-dimensional Surface Quasi-Geostrophic (SQG) equation:
\begin{equation}\label{eq:SQG}
\begin{cases}
\begin{aligned}
\partial_t u + \DIV (\bv u ) +\kappa (-\Delta)^su&=0,
\quad (\bx, t) \in D \times (0, T], \\
u(\bx, 0)
&= u_0(\bx), \quad \bx \in D,
\end{aligned}
\end{cases}
\end{equation}
where $u_0(\bx)$ is a given initial condition and periodic boundary conditions are imposed in both directions.

The SQG model is widely used in geophysical fluid dynamics, in particular for atmospheric and oceanographic simulations. In \eqref{eq:SQG} $u$ denotes the surface temperature and $\kappa \ge 0$ models the so-called Ekman pumping effect. The velocity field is defined via a fractional derivative of the solution $u$:
\begin{equation}\label{eq:neg:s}
  \bv = \nabla^\perp(-\Delta)^{-\frac12} u,
\end{equation}
where for any scalar function $w:\polR^2\to \polR$ the orthogonal gradient operator is defined as
$
\nabla^\perp w :=
\big(  -\p_{x_2}w,  \p_{x_1}w \big)^\top
$.
The negative fractional power can be computed with an approach similar to that of Section~\ref{Sec:Laps}. A detailed discussion can be found in \citep[Sec.~3.2.1]{Bonito_2021}.

\begin{figure}[t]
  \centering{
    \begin{subfigure}[t]{0.31\textwidth}
      \includegraphics[width=\textwidth]{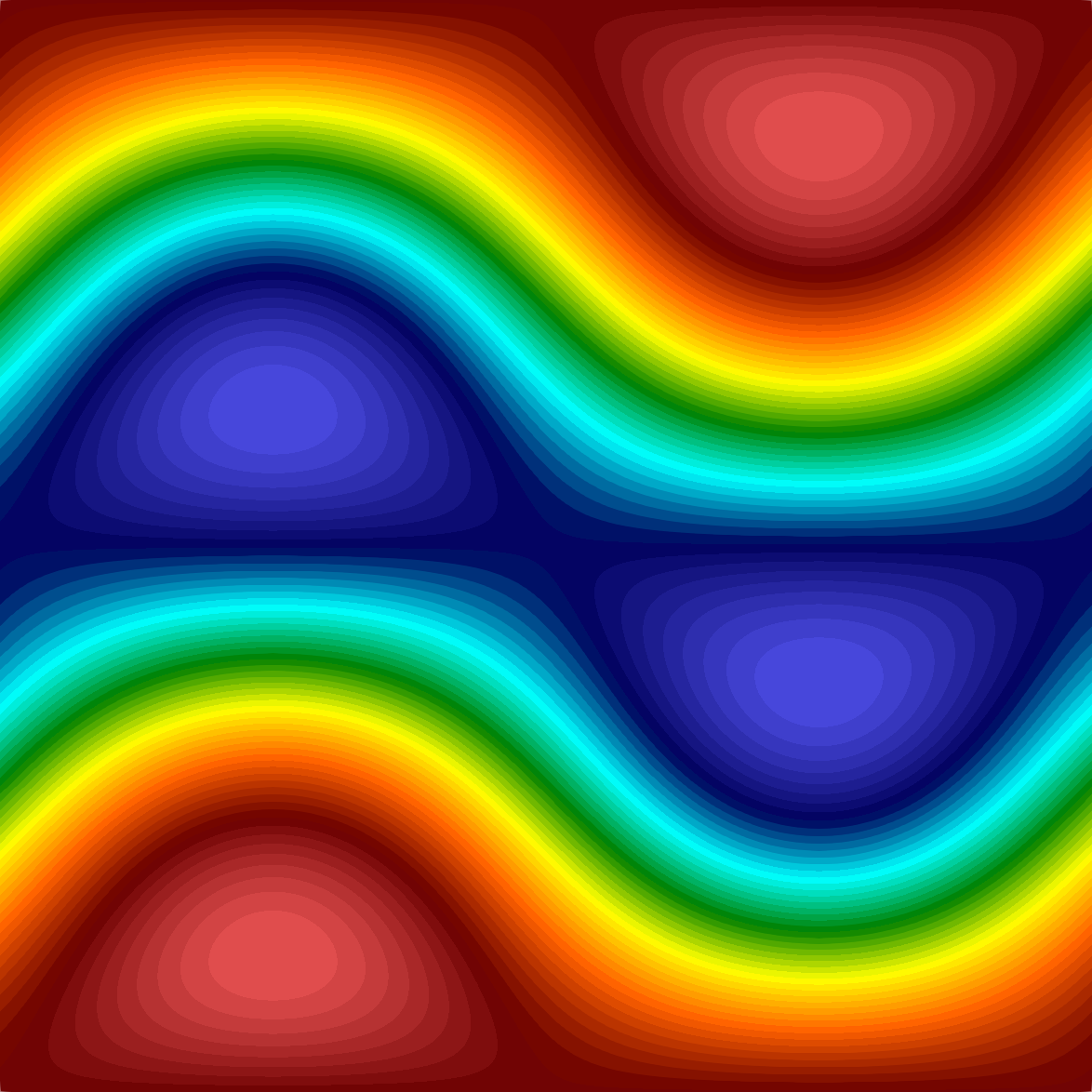}
      \caption{$t=0$}
    \end{subfigure}
    \begin{subfigure}[t]{0.31\textwidth}
      \includegraphics[width=\textwidth]{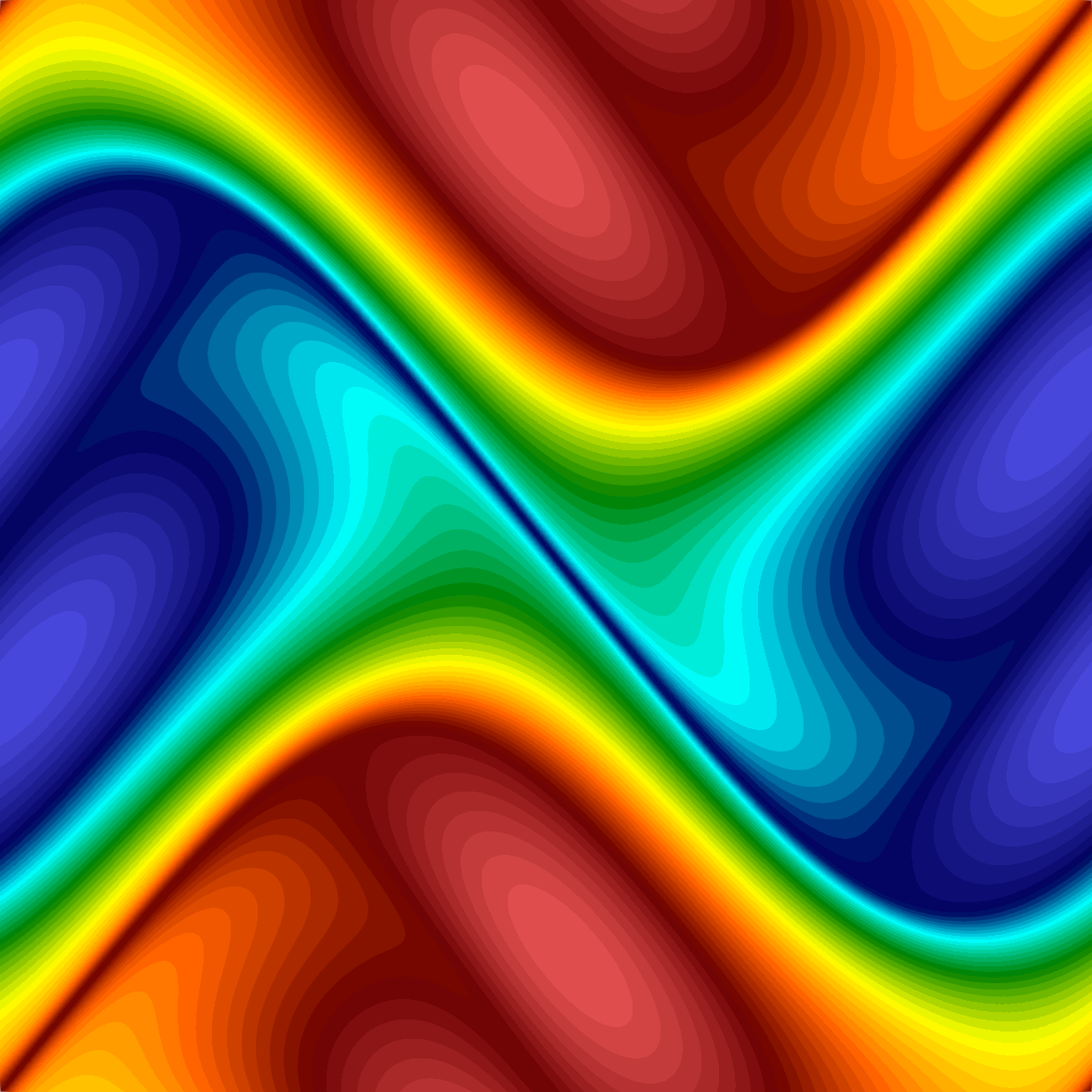}
      \caption{$t=6$}
    \end{subfigure}
    \begin{subfigure}[t]{0.31\textwidth}
      \includegraphics[width=\textwidth]{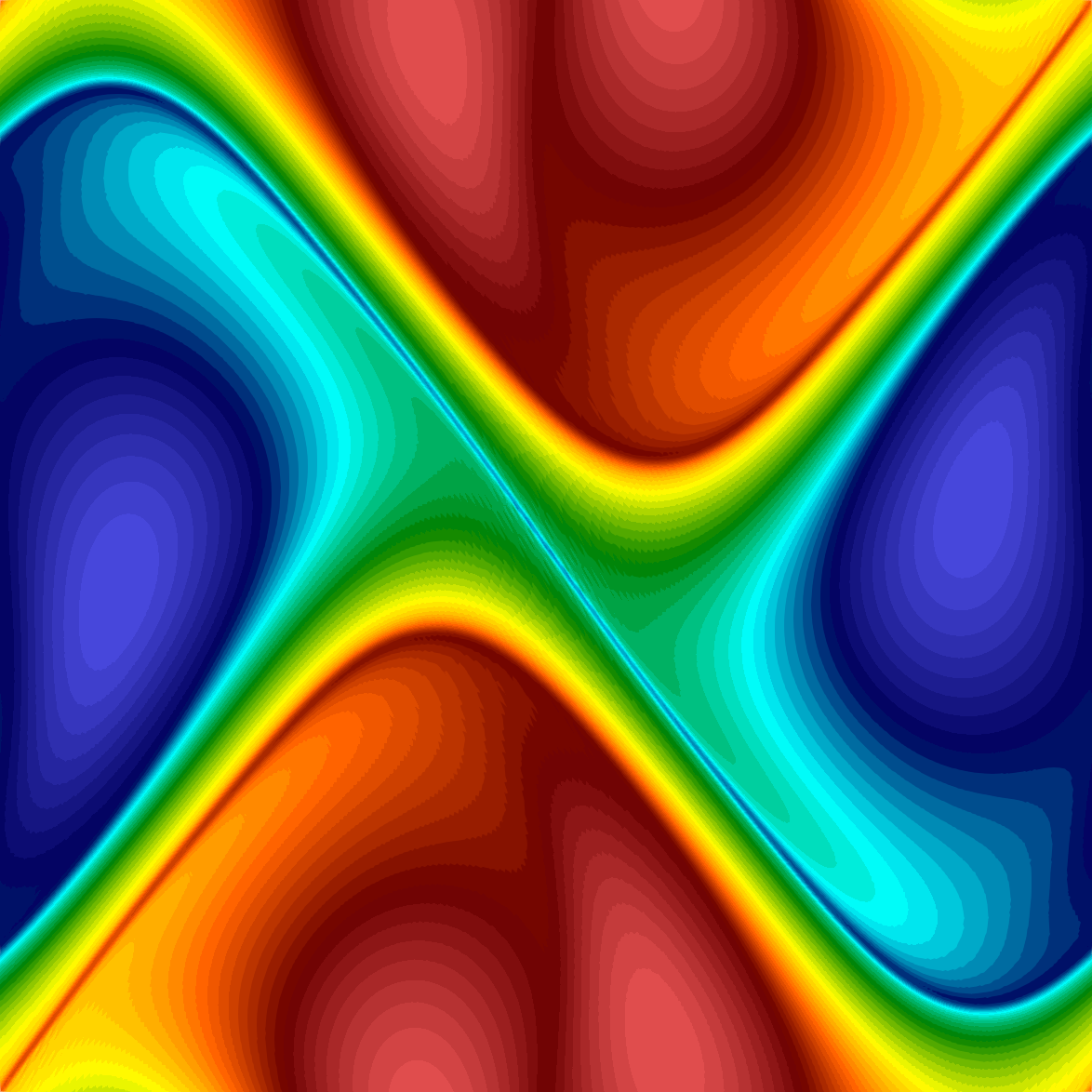}
      \caption{$t=8$}
    \end{subfigure}
    \begin{subfigure}[t]{0.0405\textwidth}
      \includegraphics[width=\textwidth]{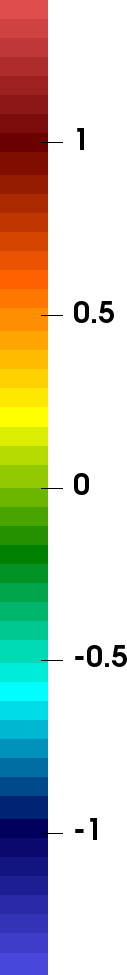}
    \end{subfigure}
  }
  \\
  \centering{
    \begin{subfigure}[t]{0.31\textwidth}
      \includegraphics[width=\textwidth]{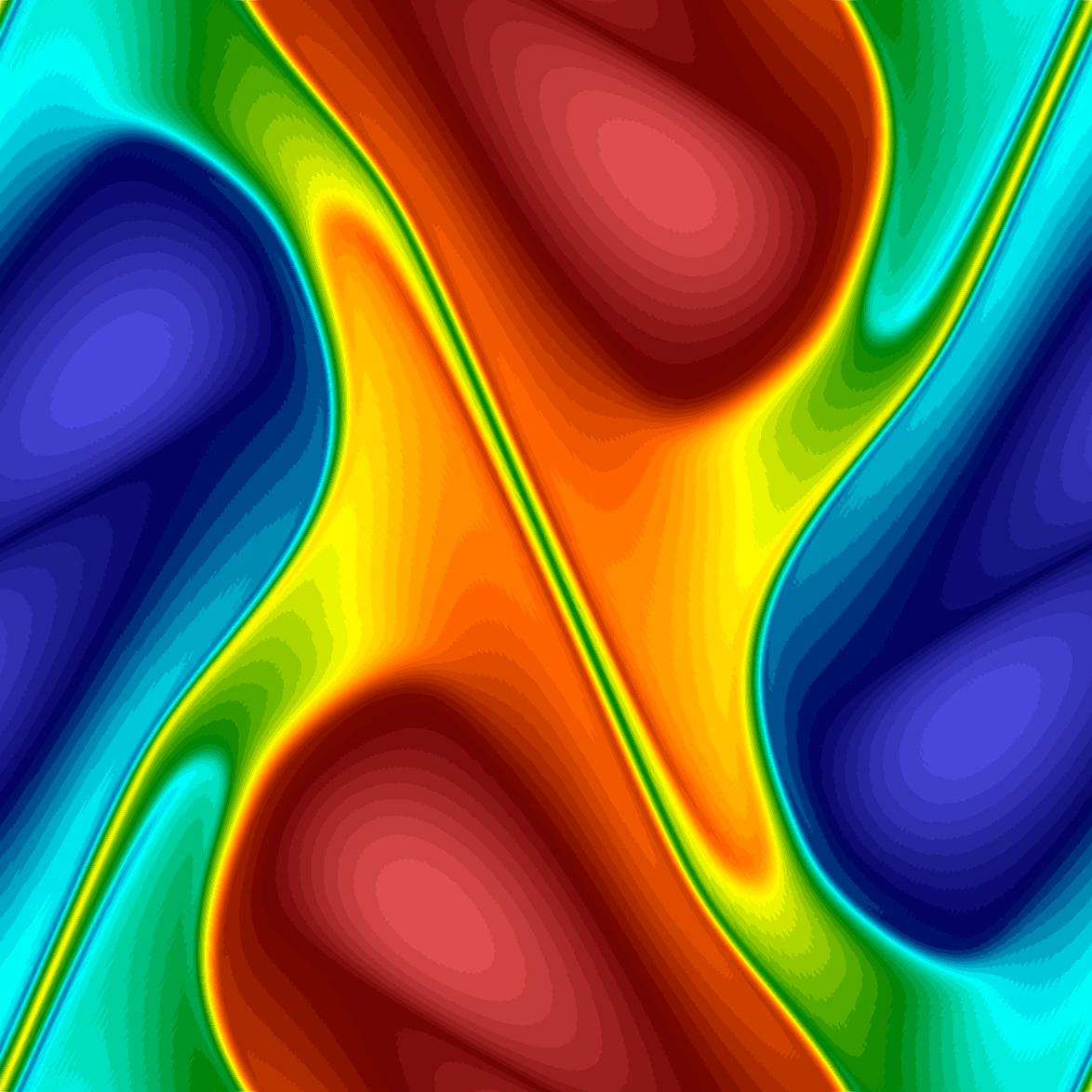}
      \caption{$t=12$}
    \end{subfigure}
    \begin{subfigure}[t]{0.31\textwidth}
      \includegraphics[width=\textwidth]{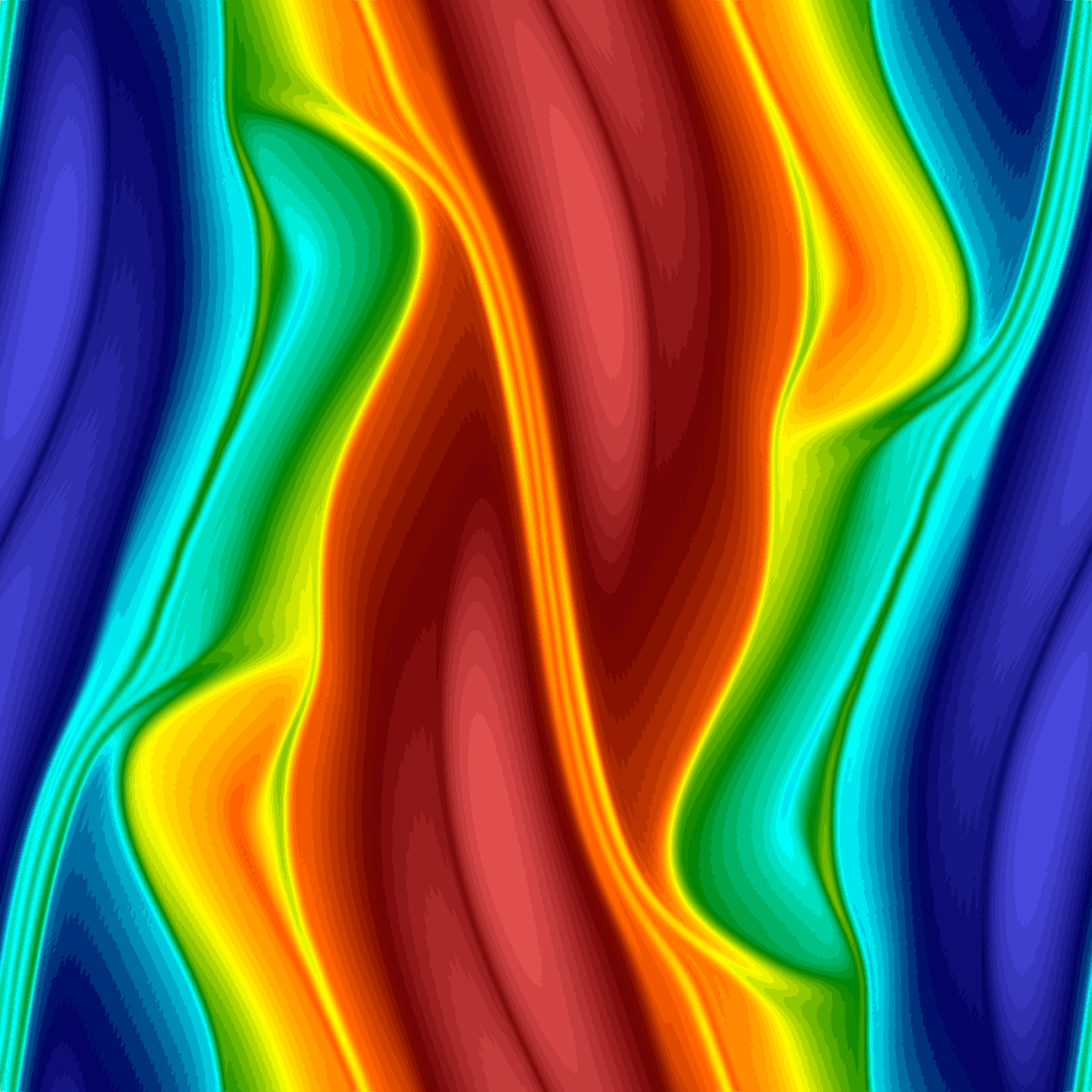}
      \caption{$t=16$}
    \end{subfigure}
    \begin{subfigure}[t]{0.31\textwidth}
      \includegraphics[width=\textwidth]{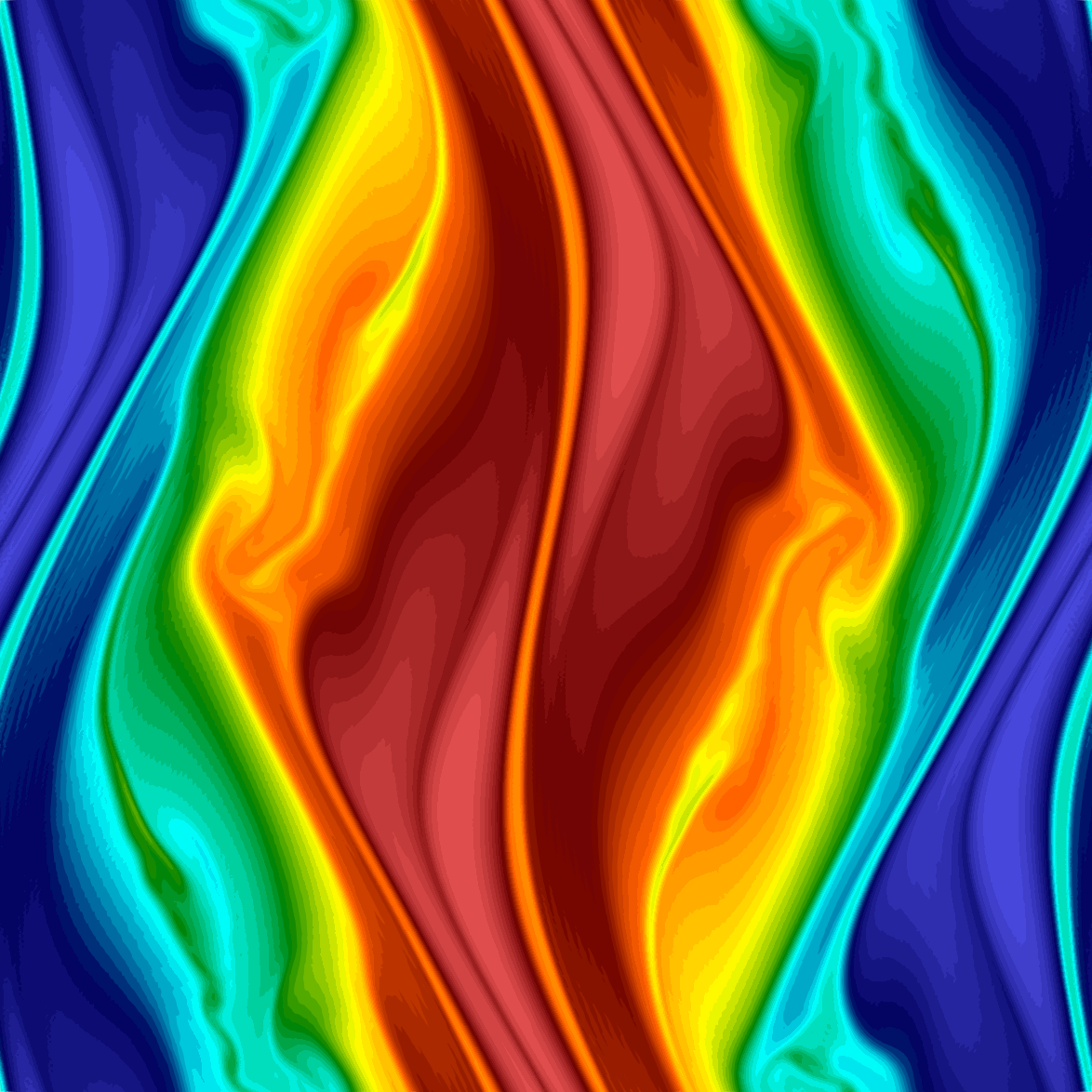}
      \caption{$t=20$}
    \end{subfigure}
    \hspace{0.175in}
  }
  \caption{Snapshots of the solution of the viscous SQG with $\kappa = 0.001$. The mesh consists of $148318$ $\polP_1$ nodes. Sharp layers are developing around time $t=8$.}
 \label{fig:sqg:k=0.001}
\end{figure}

In this final numerical demonstration, we set $D = [0, 2\pi]^2$ and use the initial condition $u_0(\bx) = \sin(x_1) \sin(x_2) + \cos(x_2)$. This initial condition is smooth, as shown in Fig.~\ref{fig:sqg:k=0.001}(a). This benchmark was previously used in \citep{Constantin_1994, Ohkitani_1997} to investigate finite-time singularity formation in the solution of the SQG equation.

Temporal snapshots of the numerical solution obtained using convex limiting with $\kappa = 0.001$ on a mesh consisting of 148318 vertices are presented in Fig.~\ref{fig:sqg:k=0.001}. We observe the development of a sharp layer around time $t \approx 8$. The algorithm successfully captures the formation of sharp layers and small-scale instabilities with mesh resolution comparable to those in the cited references, while also preserving the maximum principle.

\begin{figure}[t]
  \centering{
    \begin{subfigure}[t]{0.31\textwidth}
      \includegraphics[width=\textwidth]{FIGS/0.001/u_148318_t_8.png}
      \caption{$\kappa = 0.001$}
    \end{subfigure}
    \begin{subfigure}[t]{0.31\textwidth}
      \includegraphics[width=\textwidth]{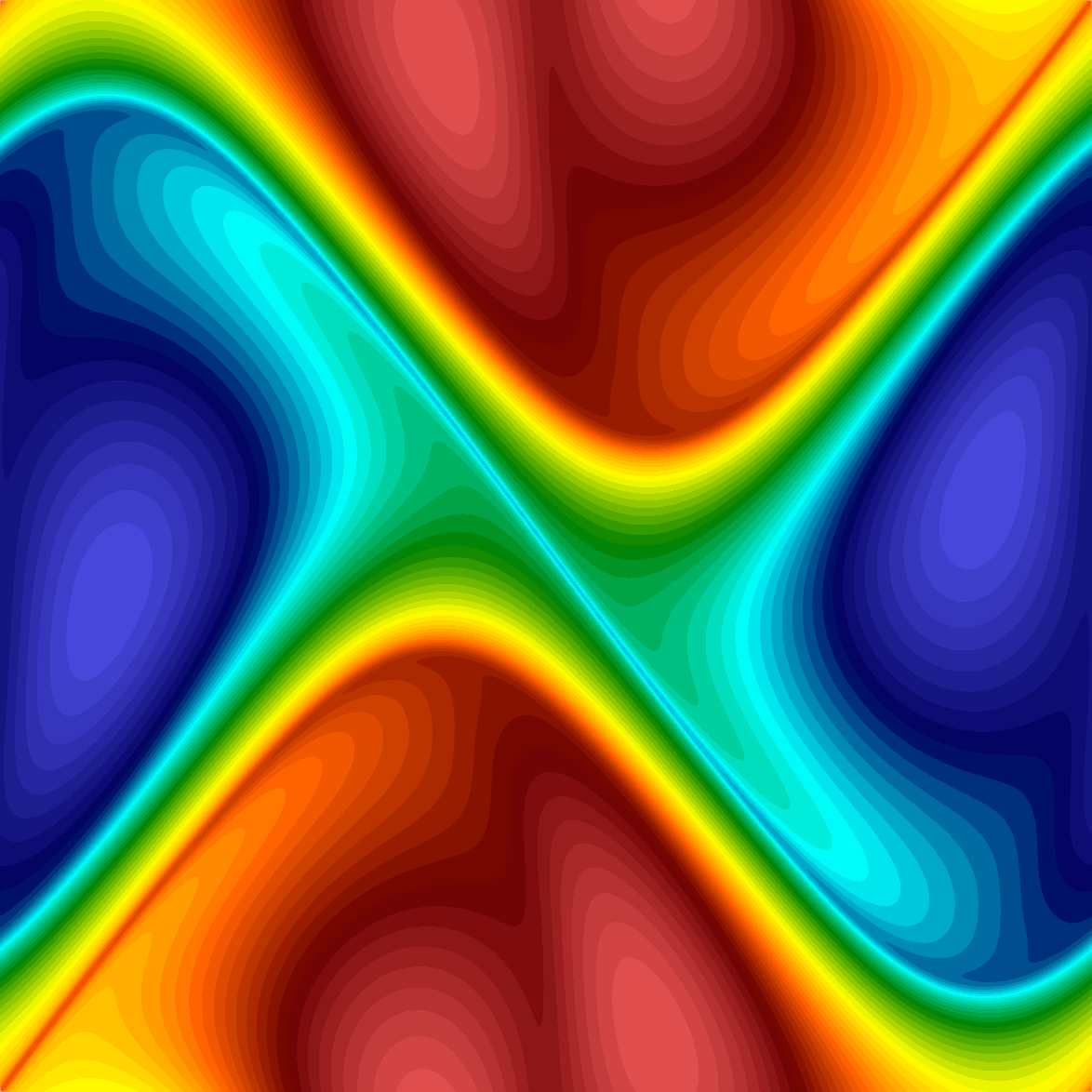}
      \caption{$\kappa = 0.01$}
    \end{subfigure}
    \begin{subfigure}[t]{0.31\textwidth}
      \includegraphics[width=\textwidth]{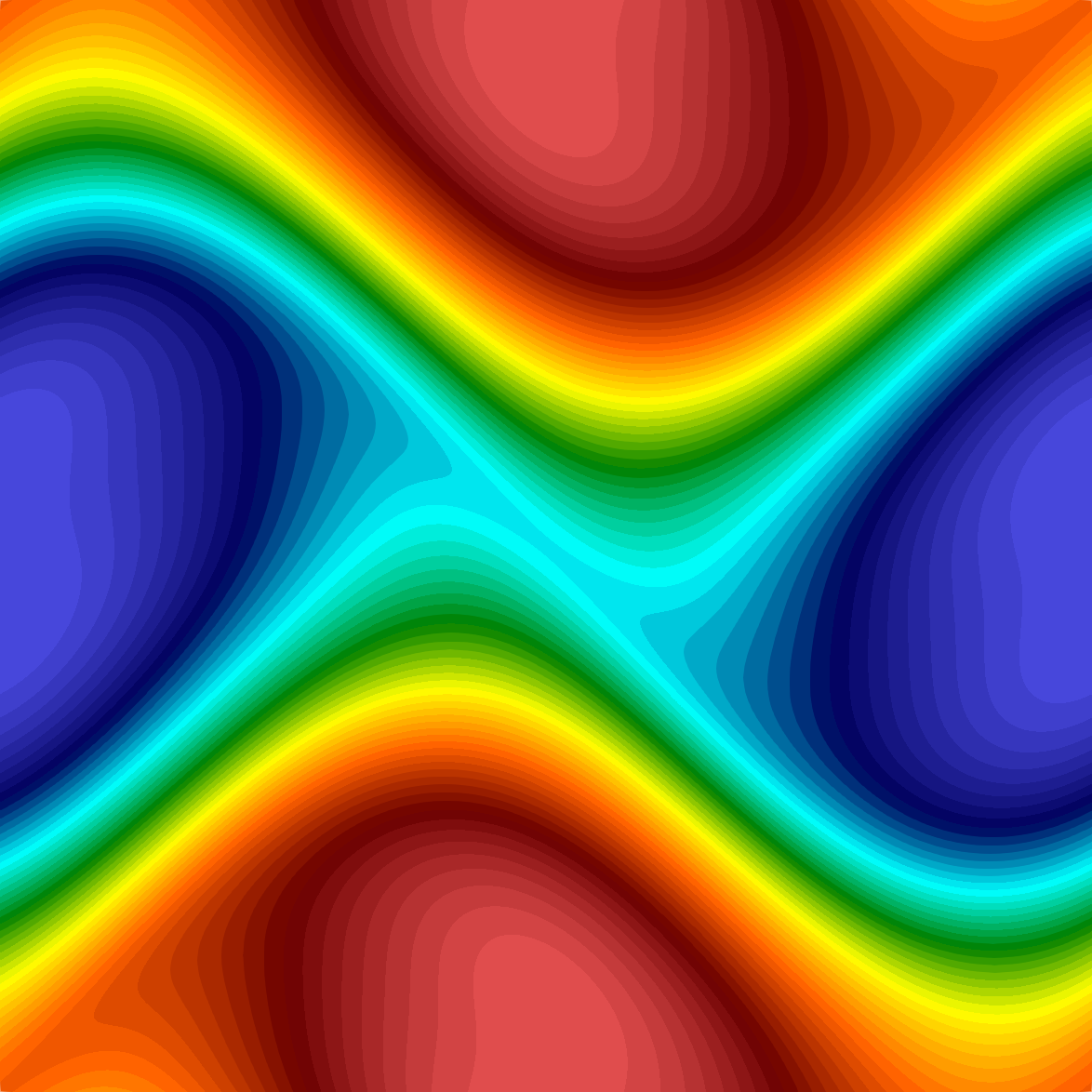}
      \caption{$\kappa = 0.1$}
    \end{subfigure}
    \begin{subfigure}[t]{0.0405\textwidth}
      \includegraphics[width=\textwidth]{FIGS/t=00_bar.png}
    \end{subfigure}
  }
  \\
  \centering{
    \begin{subfigure}[t]{0.31\textwidth}
      \includegraphics[width=\textwidth]{FIGS/0.001/u_148318_t_20.png}
      \caption{$\kappa = 0.001$}
    \end{subfigure}
    \begin{subfigure}[t]{0.31\textwidth}
      \includegraphics[width=\textwidth]{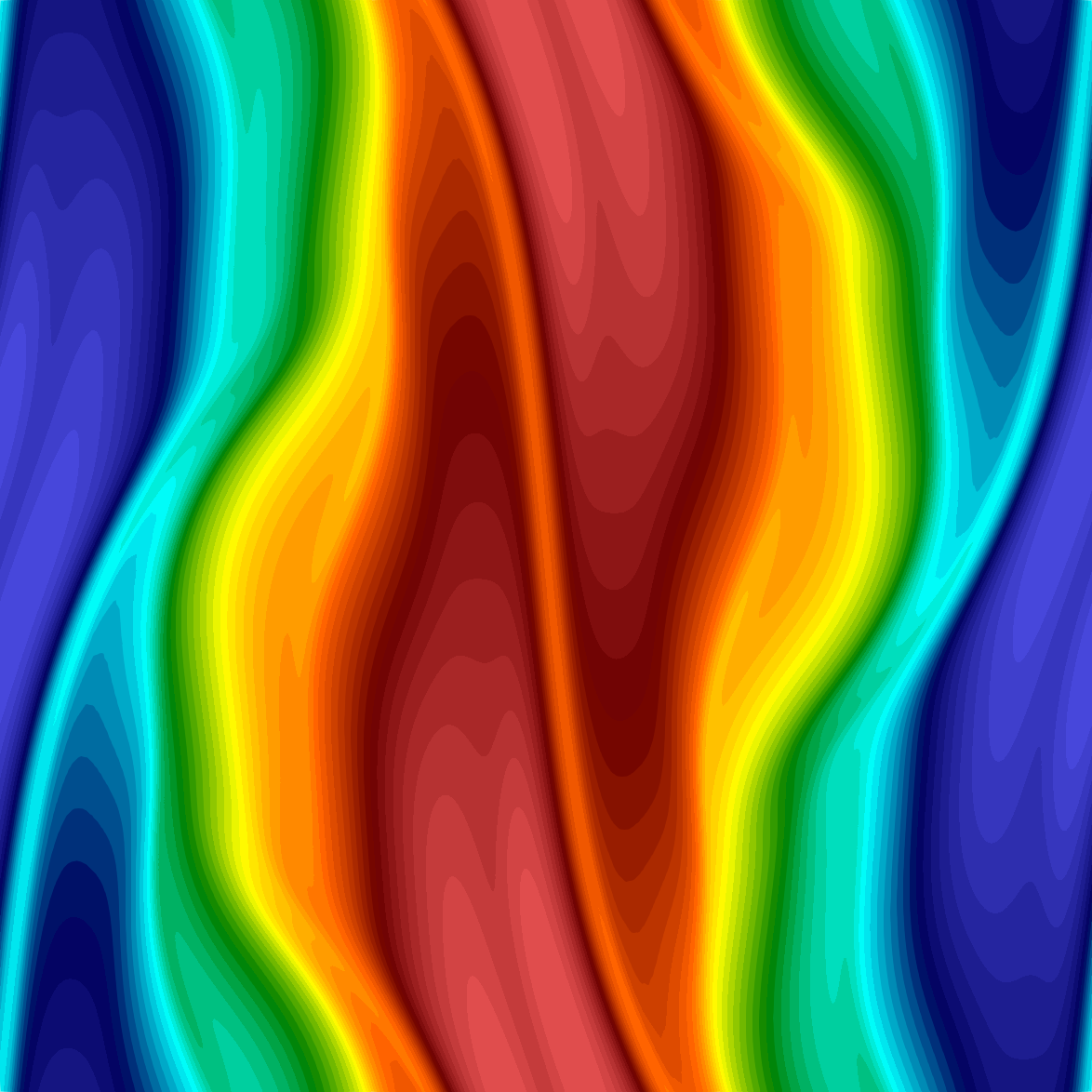}
      \caption{$\kappa = 0.01$}
    \end{subfigure}
    \begin{subfigure}[t]{0.31\textwidth}
      \includegraphics[width=\textwidth]{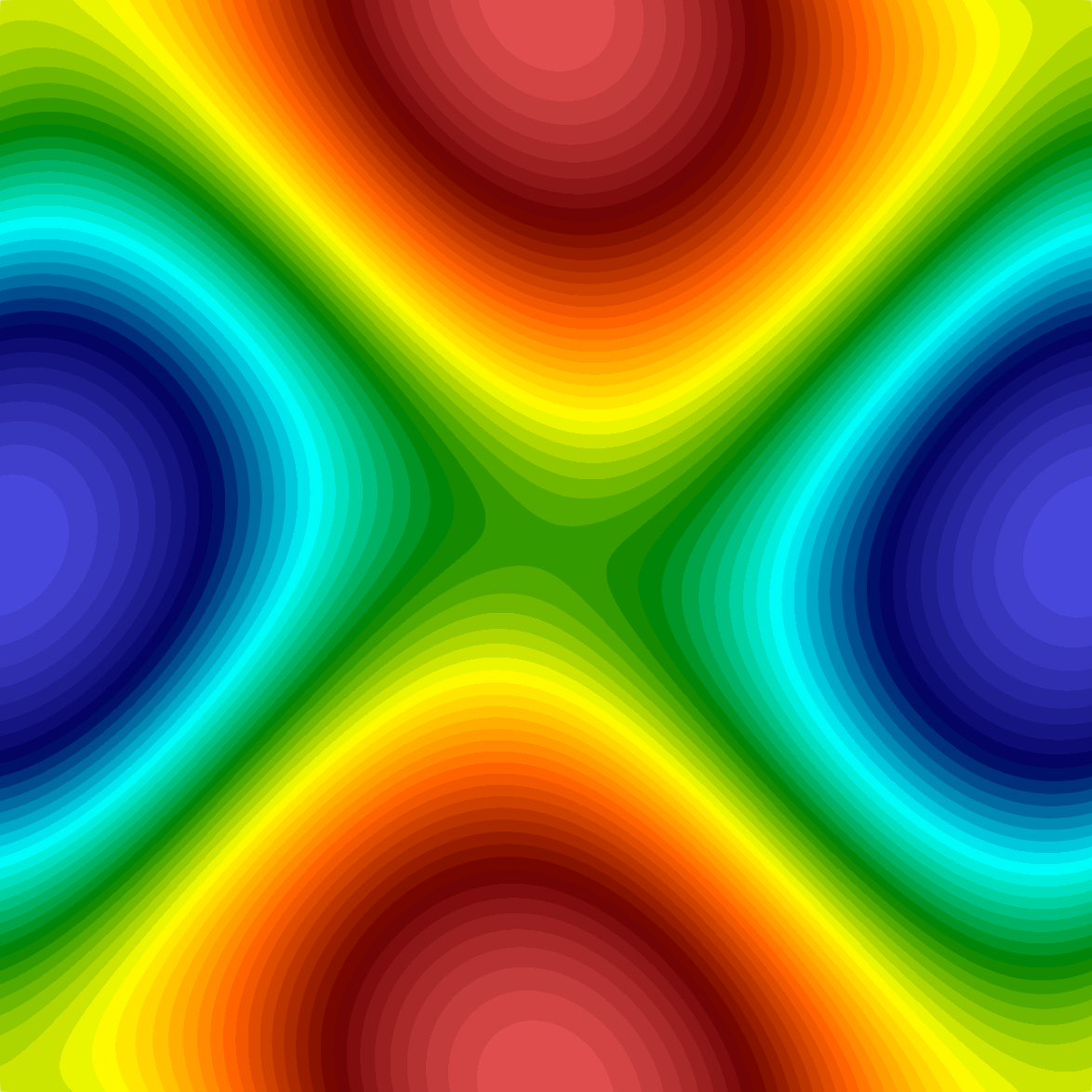}
      \caption{$\kappa = 0.1$}
    \end{subfigure}
    \hspace{0.175in}
  }
  \caption{Snapshots of the solution to the viscous SQG equation for various values of $\kappa$. The mesh consists of 148318 $\polP_1$ nodes. The top row shows the solution at time $t = 8$, while the bottom row depicts the solution at the final time $t = 20$.
}
 \label{fig:sqg}
\end{figure}

\begin{figure}[t]
  \centering{
      \includegraphics[width=0.5\textwidth]{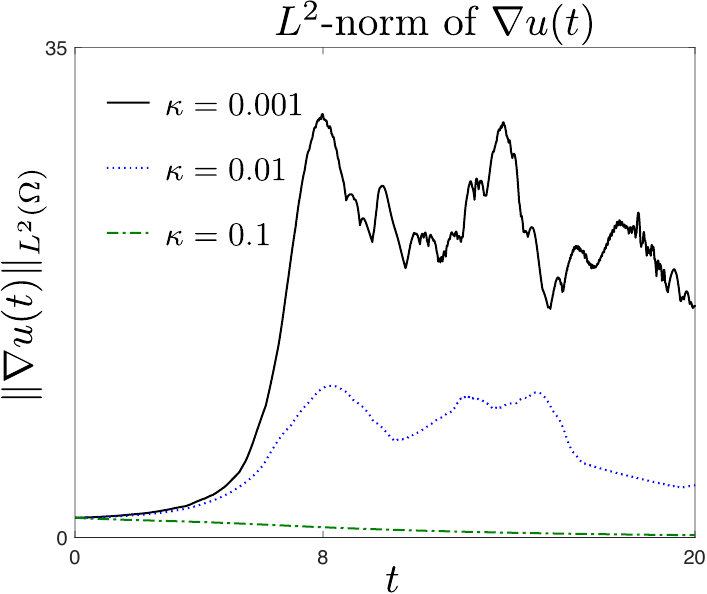}
  }
  \caption{Viscous SQG with sharp layer: Evolution of $\|\nabla u(t)\|_{L^2{(D)}}$ for different viscosity coefficients on the mesh consisting of 148318 $\polP_1$ nodes. }
  \label{fig:sqg:norm}
\end{figure}

In Fig.~\ref{fig:sqg}, we present solutions for different Ekman pumping effects. The time evolution of $\|\nabla u(t)\|_{L^2(D)}$ for various values of $\kappa$ is shown in Fig.~\ref{fig:sqg:norm}. As reported in the literature \citep{Constantin_1994,Ohkitani_1997,Bonito_2021}, we also observe a sharp increase in the gradient norm around times $t = 7$ and $t = 8$ for $\kappa = 0.001$ and $\kappa = 0.01$. For more diffusive cases, no sharp layer forms, and the gradient norm decays over time.

\FloatBarrier
\section{Concluding remarks}\label{sec:conclusions}
We have proposed a fully discrete, explicit, continuous finite element approximation of scalar conservation equations with a fractional Laplacian that preserves the maximum principle. The fractional Laplacian is approximated by a sinc quadrature of a Balakrishnan-type integral representation, which reduces its evaluation to a family of independent reaction--diffusion problems. On meshes satisfying the acute-angle condition, the discrete resolvents are monotone, and the low-order scheme with graph viscosity satisfies a discrete maximum principle under an explicit time-step restriction that accounts for both transport and fractional diffusion. A second-order approximation, which is mass conservative in the periodic setting, is then obtained by convex limiting between this low-order scheme and a high-order scheme with entropy-based viscosity. The numerical experiments confirm second-order convergence on unstructured meshes and show that the method captures the formation of sharp fronts in the SQG model while preserving the maximum principle.

The time-step restriction associated with the fractional term scales like $h^{2s}$, independently of the number of quadrature nodes; it is compatible with the hyperbolic CFL condition for $s\le\frac12$ but becomes restrictive for $s>\frac12$. Treating the fractional term implicitly while retaining the maximum principle, and relaxing the acute-angle condition, are natural directions for future work.

\section*{Declarations}
\paragraph{Funding.} This work was supported in part by NSF grant DMS-2409807 (AB) and by the Swedish Research Council (VR) grant 2021-04620 (MN).

\paragraph{Competing interests.} The authors have no relevant financial or non-financial interests to disclose.

\paragraph{Data availability.} The data supporting the findings of this study (convergence histories and figures) can be reproduced and are available from the corresponding author upon reasonable request.

\paragraph{Code availability.} The code used to produce the numerical results is available from the corresponding author upon reasonable request.

\bibliographystyle{abbrvnat}

\bibliography{ref_jsc}

\begin{thebibliography}{36}
\providecommand{\natexlab}[1]{#1}
\providecommand{\url}[1]{\texttt{#1}}
\expandafter\ifx\csname urlstyle\endcsname\relax
  \providecommand{\doi}[1]{doi: #1}\else
  \providecommand{\doi}{doi: \begingroup \urlstyle{rm}\Url}\fi

\bibitem[Bonito and Lei(2022)]{bonito2021approximation}
A.~Bonito and W.~Lei.
\newblock Approximation of the spectral fractional powers of the
  {L}aplace--{B}eltrami operator.
\newblock \emph{Numer. Math. Theory Methods Appl.}, 15\penalty0 (4):\penalty0
  1193--1218, 2022.
\newblock \doi{10.4208/nmtma.OA-2022-0005s}.

\bibitem[Bonito and Nazarov(2021)]{Bonito_2021}
A.~Bonito and M.~Nazarov.
\newblock Numerical simulations of surface quasi-geostrophic flows on periodic
  domains.
\newblock \emph{SIAM J. Sci. Comput.}, 43\penalty0 (2):\penalty0 B405--B430,
  2021.
\newblock ISSN 1064-8275,1095-7197.
\newblock \doi{10.1137/20M1342616}.
\newblock URL \url{https://doi.org/10.1137/20M1342616}.

\bibitem[Bonito and Pasciak(2015)]{bonito2015numerical}
A.~Bonito and J.~Pasciak.
\newblock Numerical approximation of fractional powers of elliptic operators.
\newblock \emph{Mathematics of Computation}, 84\penalty0 (295):\penalty0
  2083--2110, 2015.

\bibitem[Bonito and Pasciak(2017)]{bonito2017numerical}
A.~Bonito and J.~E. Pasciak.
\newblock Numerical approximation of fractional powers of regularly accretive
  operators.
\newblock \emph{IMA Journal of Numerical Analysis}, 37\penalty0 (3):\penalty0
  1245--1273, 2017.

\bibitem[Bonito and Wei(2020)]{bonito2020electroconvection}
A.~Bonito and P.~Wei.
\newblock Electroconvection of thin liquid crystals: Model reduction and
  numerical simulations.
\newblock \emph{Journal of Computational Physics}, 405:\penalty0 109140, 2020.

\bibitem[Bonito et~al.(2018)Bonito, Borthagaray, Nochetto, Ot{\'a}rola, and
  Salgado]{bonito2018numerical}
A.~Bonito, J.~P. Borthagaray, R.~H. Nochetto, E.~Ot{\'a}rola, and A.~J.
  Salgado.
\newblock Numerical methods for fractional diffusion.
\newblock \emph{Computing and Visualization in Science}, 19\penalty0
  (5):\penalty0 19--46, 2018.

\bibitem[Bonito et~al.(2019{\natexlab{a}})Bonito, Lei, and
  Pasciak]{bonito2017sinc}
A.~Bonito, W.~Lei, and J.~E. Pasciak.
\newblock On sinc quadrature approximations of fractional powers of regularly
  accretive operators.
\newblock \emph{J. Numer. Math.}, 27\penalty0 (2):\penalty0 57--68,
  2019{\natexlab{a}}.
\newblock \doi{10.1515/jnma-2017-0116}.

\bibitem[Bonito et~al.(2019{\natexlab{b}})Bonito, Lei, and
  Pasciak]{bonito2019numerical}
A.~Bonito, W.~Lei, and J.~E. Pasciak.
\newblock Numerical approximation of the integral fractional {L}aplacian.
\newblock \emph{Numerische Mathematik}, 142\penalty0 (2):\penalty0 235--278,
  2019{\natexlab{b}}.

\bibitem[Brandts et~al.(2009)Brandts, Korotov, K{\v{r}}{\'{\i}}{\v{z}}ek, and
  {\v{S}}olc]{BrKKS:09}
J.~Brandts, S.~Korotov, M.~K{\v{r}}{\'{\i}}{\v{z}}ek, and J.~{\v{S}}olc.
\newblock On nonobtuse simplicial partitions.
\newblock \emph{SIAM Rev.}, 51\penalty0 (2):\penalty0 317--335, 2009.

\bibitem[Brandts et~al.(2008)Brandts, Korotov, and
  K{\v{r}}{\'{\i}}{\v{z}}ek]{BrKoK:08}
J.~H. Brandts, S.~Korotov, and M.~K{\v{r}}{\'{\i}}{\v{z}}ek.
\newblock The discrete maximum principle for linear simplicial finite element
  approximations of a reaction-diffusion problem.
\newblock \emph{Linear Algebra Appl.}, 429\penalty0 (10):\penalty0 2344--2357,
  2008.

\bibitem[Clavin(2002)]{Clavin_2002}
P.~Clavin.
\newblock Instabilities and nonlinear patterns of overdriven detonations in
  gases.
\newblock In H.~Berestycki and Y.~Pomeau, editors, \emph{Nonlinear {PDE}'s in
  Condensed Matter and Reactive Flows}, pages 49--97. Springer Netherlands,
  Dordrecht, 2002.
\newblock ISBN 978-94-010-0307-0.
\newblock \doi{10.1007/978-94-010-0307-0_3}.
\newblock URL \url{https://doi.org/10.1007/978-94-010-0307-0_3}.

\bibitem[Constantin et~al.(1994)Constantin, Majda, and Tabak]{Constantin_1994}
P.~Constantin, A.~J. Majda, and E.~Tabak.
\newblock Formation of strong fronts in the {$2$}-{D} quasigeostrophic thermal
  active scalar.
\newblock \emph{Nonlinearity}, 7\penalty0 (6):\penalty0 1495--1533, 1994.
\newblock ISSN 0951-7715,1361-6544.
\newblock \doi{10.1088/0951-7715/7/6/001}.
\newblock URL \url{https://doi.org/10.1088/0951-7715/7/6/001}.

\bibitem[Constantin et~al.(2012)Constantin, Lai, Sharma, Tseng, and
  Wu]{Constantin_2012}
P.~Constantin, M.-C. Lai, R.~Sharma, Y.-H. Tseng, and J.~Wu.
\newblock New numerical results for the surface quasi-geostrophic equation.
\newblock \emph{J. Sci. Comput.}, 50\penalty0 (1):\penalty0 1--28, 2012.
\newblock ISSN 0885-7474,1573-7691.
\newblock \doi{10.1007/s10915-011-9471-9}.
\newblock URL \url{https://doi.org/10.1007/s10915-011-9471-9}.

\bibitem[del Castillo-Negrete(2010)]{del_Castillo_2010}
D.~del Castillo-Negrete.
\newblock Non-diffusive, non-local transport in fluids and plasmas.
\newblock \emph{Nonlinear Processes in Geophysics}, 17\penalty0 (6):\penalty0
  795--807, 2010.
\newblock \doi{10.5194/npg-17-795-2010}.
\newblock URL \url{https://npg.copernicus.org/articles/17/795/2010/}.

\bibitem[D'Elia et~al.(2020)D'Elia, Du, Glusa, Gunzburger, Tian, and
  Zhou]{d2020numerical}
M.~D'Elia, Q.~Du, C.~Glusa, M.~Gunzburger, X.~Tian, and Z.~Zhou.
\newblock Numerical methods for nonlocal and fractional models.
\newblock \emph{Acta Numerica}, 29:\penalty0 1--124, 2020.

\bibitem[Ern and Guermond(2021)]{Ern_Guermond_2021}
A.~Ern and J.-L. Guermond.
\newblock \emph{Finite elements {III}---first-order and time-dependent {PDE}s},
  volume~74 of \emph{Texts in Applied Mathematics}.
\newblock Springer, Cham, 2021.
\newblock ISBN 978-3-030-57347-8; 978-3-030-57348-5.
\newblock \doi{10.1007/978-3-030-57348-5}.
\newblock URL \url{https://doi.org/10.1007/978-3-030-57348-5}.

\bibitem[Ferracina and Spijker(2005)]{Ferracina_Spijker_2005}
L.~Ferracina and M.~N. Spijker.
\newblock An extension and analysis of the {S}hu-{O}sher representation of
  {R}unge-{K}utta methods.
\newblock \emph{Math. Comp.}, 74\penalty0 (249):\penalty0 201--219, 2005.

\bibitem[Geuzaine and Remacle(2009)]{gmsh_2009}
C.~Geuzaine and J.-F. Remacle.
\newblock Gmsh: {A} 3-{D} finite element mesh generator with built-in pre- and
  post-processing facilities.
\newblock \emph{Internat. J. Numer. Methods Engrg.}, 79\penalty0 (11):\penalty0
  1309--1331, 2009.
\newblock ISSN 0029-5981,1097-0207.
\newblock \doi{10.1002/nme.2579}.
\newblock URL \url{https://doi.org/10.1002/nme.2579}.

\bibitem[Gottlieb et~al.(2009)Gottlieb, Ketcheson, and
  Shu]{Gottlieb_Ketcheson_Shu_2009}
S.~Gottlieb, D.~I. Ketcheson, and C.-W. Shu.
\newblock High order strong stability preserving time discretizations.
\newblock \emph{J. Sci. Comput.}, 38\penalty0 (3):\penalty0 251--289, 2009.

\bibitem[Guermond and Nazarov(2014)]{Guermond_Nazarov_2014}
J.-L. Guermond and M.~Nazarov.
\newblock A maximum-principle preserving {$C^0$} finite element method for
  scalar conservation equations.
\newblock \emph{Comput. Methods Appl. Mech. Engrg.}, 272:\penalty0 198--213,
  2014.
\newblock ISSN 0045-7825,1879-2138.
\newblock \doi{10.1016/j.cma.2013.12.015}.
\newblock URL \url{https://doi.org/10.1016/j.cma.2013.12.015}.

\bibitem[Guermond and Popov(2016)]{Guermond_Popov_2016}
J.-L. Guermond and B.~Popov.
\newblock Invariant domains and first-order continuous finite element
  approximation for hyperbolic systems.
\newblock \emph{SIAM J. Numer. Anal.}, 54\penalty0 (4):\penalty0 2466--2489,
  2016.

\bibitem[Guermond and Popov(2017)]{guermond_popov_second_order_2018}
J.-L. Guermond and B.~Popov.
\newblock Invariant {D}omains and {S}econd-{O}rder {C}ontinuous {F}inite
  {E}lement {A}pproximation for {S}calar {C}onservation {E}quations.
\newblock \emph{SIAM J. Numer. Anal.}, 55\penalty0 (6):\penalty0 3120--3146,
  2017.

\bibitem[Guermond et~al.(2011)Guermond, Pasquetti, and
  Popov]{guermond2011entropy}
J.-L. Guermond, R.~Pasquetti, and B.~Popov.
\newblock Entropy viscosity method for nonlinear conservation laws.
\newblock \emph{Journal of Computational Physics}, 230\penalty0 (11):\penalty0
  4248--4267, 2011.

\bibitem[Guermond et~al.(2018)Guermond, Nazarov, Popov, and
  Tomas]{Guermond_Nazarov_Popov_Tomas_2018}
J.-L. Guermond, M.~Nazarov, B.~Popov, and I.~Tomas.
\newblock Second-order invariant domain preserving approximation of the {E}uler
  equations using convex limiting.
\newblock \emph{SIAM J. Sci. Comput.}, 40\penalty0 (5):\penalty0 A3211--A3239,
  2018.

\bibitem[Higueras(2005)]{Higueras_2005}
I.~Higueras.
\newblock Representations of {R}unge-{K}utta methods and strong stability
  preserving methods.
\newblock \emph{SIAM J. Numer. Anal.}, 43\penalty0 (3):\penalty0 924--948,
  2005.

\bibitem[Huang and Oberman(2014)]{Huang_2014}
Y.~Huang and A.~Oberman.
\newblock Numerical methods for the fractional {L}aplacian: a finite
  difference--quadrature approach.
\newblock \emph{SIAM J. Numer. Anal.}, 52\penalty0 (6):\penalty0 3056--3084,
  2014.
\newblock ISSN 0036-1429,1095-7170.
\newblock \doi{10.1137/140954040}.
\newblock URL \url{https://doi.org/10.1137/140954040}.

\bibitem[Logg et~al.(2012)Logg, Mardal, and Wells]{fenics}
A.~Logg, K.-A. Mardal, and G.~Wells.
\newblock \emph{Automated solution of differential equations by the finite
  element method: The FEniCS book}.
\newblock Springer Science \& Business Media, 2012.

\bibitem[Lunardi(2018)]{lunardi2018interpolation}
A.~Lunardi.
\newblock \emph{Interpolation Theory}, volume~16 of \emph{Lecture Notes. Scuola
  Normale Superiore di Pisa (New Series)}.
\newblock Edizioni della Normale, Pisa, third edition, 2018.

\bibitem[Lund and Bowers(1992)]{KBSinc92}
J.~Lund and K.~L. Bowers.
\newblock \emph{Sinc methods for quadrature and differential equations}.
\newblock Society for Industrial and Applied Mathematics (SIAM), Philadelphia,
  PA, 1992.
\newblock ISBN 0-89871-298-X.
\newblock \doi{10.1137/1.9781611971637}.
\newblock URL \url{https://doi.org/10.1137/1.9781611971637}.

\bibitem[Nochetto et~al.(2015)Nochetto, Ot{\'a}rola, and
  Salgado]{nochetto2015pde}
R.~H. Nochetto, E.~Ot{\'a}rola, and A.~J. Salgado.
\newblock A {PDE} approach to fractional diffusion in general domains: a priori
  error analysis.
\newblock \emph{Found. Comput. Math.}, 15\penalty0 (3):\penalty0 733--791,
  2015.
\newblock \doi{10.1007/s10208-014-9208-x}.

\bibitem[Ohkitani and Yamada(1997)]{Ohkitani_1997}
K.~Ohkitani and M.~Yamada.
\newblock Inviscid and inviscid-limit behavior of a surface quasigeostrophic
  flow.
\newblock \emph{Phys. Fluids}, 9\penalty0 (4):\penalty0 876--882, 1997.
\newblock ISSN 1070-6631,1089-7666.
\newblock \doi{10.1063/1.869184}.
\newblock URL \url{https://doi.org/10.1063/1.869184}.

\bibitem[Shu and Osher(1988)]{Shu_Osher1988}
C.-W. Shu and S.~Osher.
\newblock Efficient implementation of essentially non-oscillatory
  shock-capturing schemes.
\newblock \emph{J. Comput. Phys.}, 77\penalty0 (2):\penalty0 439--471, 1988.

\bibitem[Song and Karniadakis(2017)]{Song_2017}
F.~Song and G.~E. Karniadakis.
\newblock Fractional spectral vanishing viscosity method: application to the
  quasi-geostrophic equation.
\newblock \emph{Chaos Solitons Fractals}, 102:\penalty0 327--332, 2017.
\newblock ISSN 0960-0779,1873-2887.
\newblock \doi{10.1016/j.chaos.2017.03.052}.
\newblock URL \url{https://doi.org/10.1016/j.chaos.2017.03.052}.

\bibitem[Sugimoto(1991)]{Sugimoto_1991}
N.~Sugimoto.
\newblock Burgers equation with a fractional derivative; hereditary effects on
  nonlinear acoustic waves.
\newblock \emph{J. Fluid Mech.}, 225:\penalty0 631--653, 1991.
\newblock ISSN 0022-1120,1469-7645.
\newblock \doi{10.1017/S0022112091002203}.
\newblock URL \url{https://doi.org/10.1017/S0022112091002203}.

\bibitem[Woyczy{\'n}ski(2001)]{Woyczynski_2001}
W.~A. Woyczy{\'n}ski.
\newblock L{\'e}vy processes in the physical sciences.
\newblock In O.~E. Barndorff-Nielsen, T.~Mikosch, and S.~I. Resnick, editors,
  \emph{L{\'e}vy Processes: Theory and Applications}, pages 241--266.
  Birkh{\"a}user, Boston, MA, 2001.
\newblock \doi{10.1007/978-1-4612-0197-7_11}.

\bibitem[Xu and Zikatanov(1999)]{XuZik:99}
J.~Xu and L.~Zikatanov.
\newblock A monotone finite element scheme for convection-diffusion equations.
\newblock \emph{Math. Comp.}, 68\penalty0 (228):\penalty0 1429--1446, 1999.

\end{thebibliography}

\end{document}